\documentclass[11pt]{article}

\usepackage{amsmath}
\usepackage{amssymb}
\usepackage[latin1]{inputenc}

\makeatletter
\newlength{\earraycolsep}
\def\eqnarray{\stepcounter{equation}\let\@currentlabel%
\theequation
\global\@eqnswtrue\m@th
\global\@eqcnt\z@\tabskip\@centering\let\\\@eqncr
$$\halign to\displaywidth\bgroup\@eqnsel\hskip\@centering
$\displaystyle\tabskip\z@{##}$&\global\@eqcnt\@ne
\hskip 2\earraycolsep \hfil$\displaystyle{##}$\hfil
&\global\@eqcnt\tw@ \hskip 2\earraycolsep
$\displaystyle\tabskip\z@{##}$\hfil
\tabskip\@centering&\llap{##}\tabskip\z@\cr}
\makeatother

\newcommand{\R}{\mathbb R}

\newcommand{\n}{\mathbf{n}}

\renewcommand{\AA}{\mathfrak{A}}

\newcommand{\ct}[1]{\mathbf{#1}_{\tau}}
\newcommand{\h}{{\bf{h}}}
\def\nablat{\nabla_{\tau}}
\def\Deltat{\Delta_{\tau}}
\def\Deltatt{\Delta_{\tau,t}}
\def\phiet{\tilde{\phi}}

\def\gammat{\gamma_{\tau}}
\renewcommand{\div}[1]{\mbox{\rm div}\left(#1\right)}
\newcommand{\divt}[1]{\mbox{\rm div}_{\tau}\left(#1\right)}
\newcommand{\ud}{u_d}
\newcommand{\un}{u_n}
\newcommand{\vd}{v_d}
\newcommand{\vn}{v_n}

\newcommand{\KV}{J_{KV}}

\newcommand{\sH}{{\rm H}}
\newcommand{\sL}{{\rm L}}

\newtheorem{theorem}{Theorem}
\newtheorem{proposition}{Proposition}

\newtheorem{lemma}{Lemma}
\newtheorem{definition}{Definition}
\newtheorem{remark}{Remark}

\newenvironment{proofof}[1]{\begin{trivlist}
                       \item[]\hspace{0cm}{\bf Proof of #1 }
                       \hspace{0cm} }{\hfill $\blacksquare$
                     \end{trivlist}}

\begin{document}

\title{On second order shape optimization methods for electrical impedance tomography.}
\author{L. Afraites, M. Dambrine and D. Kateb\footnote{Laboratoire de Math\'ematiques Appliqu\'ees de Compi\`egne. Universit\'e de Technologie de Compi\`egne.}}
%\date{}
\maketitle

\abstract{This paper is devoted to the analysis of a second order method for recovering the \emph{a priori} unknown shape of an inclusion $\omega$ inside a body $\Omega$ from boundary measurement. This inverse problem - known as electrical impedance tomography - has many important practical applications and hence has focussed much attention during the last years. However, to our best knowledge, no work has yet  considered a second order approach for this problem. This paper aims to fill that void: we investigate the existence of second order derivative of the state $u$ with respect to perturbations of the shape of the interface $\partial\omega$, then we choose a cost function in order to recover the geometry of $\partial \omega$ and derive the expression of the derivatives needed to implement the corresponding Newton method. We then investigate the stability of the process and explain why this inverse problem is severely ill-posed by proving the compactness of the Hessian at the global minimizer.}

%%%%%%%%%%%%%%%%%%%%%%%%%%%%%
\section{Introduction and statement of the results.}
%%%%%%%%%%%%%%%%%%%%%%%%%%%%%

Let $\Omega$ be a bounded open set with smooth boundary in $\R^2$ or $\R^3$.
Consider a $L^{\infty}$ function $\sigma$ such that there exists a real $c$
with $\sigma(x) \geq c>0$. Consider the elliptic equation
\begin{equation}
\label{equation:etat}
-\div{ \sigma(x)  \nabla u } =0 \text{ in } \Omega,
\end{equation}
with the Dirichlet boundary condition
\begin{equation}
\label{equation:potentiel:condition:limite:dirichlet}
u = f \text{ on } \partial \Omega.
\end{equation}
Define the Dirichlet-to-Neumann map as
\begin{equation*}
\Lambda_\sigma : f \mapsto \sigma \left( \partial_\n u\right)_{|\partial \Omega},
\end{equation*}
where $u$ solves \eqref{equation:etat},\eqref{equation:potentiel:condition:limite:dirichlet} and $\n$ is the outer unit normal vector to $\partial\Omega$. The inverse conductivity problem of Calder\'on is to determine $\sigma$ from $\Lambda_{\sigma}$. Electrical impedance tomography aims to form an image of the conductivity distribution $\sigma$ from the knowledge of $\Lambda_\sigma$. When $\sigma$ is smooth enough, one can reconstruct $\sigma$ from $\Lambda_{\sigma}$ (see the works of Sylvester and Uhlmann \cite{SylvesterUhlmann}, Nachmann \cite{Nachmann1,Nachmann2} and Novikov
\cite{Novikov}). When the conductivity distribution is only $\sL^{\infty}$, Astala and P\"{a}iv\"{a}rinta have recently shown in \cite{AstalaPaivarinta} that, in dimension two, the map $\Lambda_{\sigma}$ determines $\sigma \in \sL^{\infty}(\Omega)$.

We are interested in a particular case of that  problem:  when a body is inserted inside a given object with a distinct conductivity, the question of determining its shape from boundary measurement arises in many fields of modern technology. In the context of the  inverse problem of conductivity of Calder\'on, we restrict the range of admissible conductivity distributions to the family of  piecewise constant functions which take only  two distinct values $\sigma_{1},\sigma_2>0$ which are assumed to be known. The conductivity distribution is then defined by an open subset $\omega$  as
\begin{equation}
\label{definition:conductivite}
\sigma = \sigma_{1} \chi_{\Omega \setminus \omega}+ \sigma_{2} \chi_{\omega}.
\end{equation}
Here, the only unknown of the problem is $\omega$  a subdomain of $\Omega$ with a smooth boundary $\partial \omega$; its outer unit normal vector is denoted by $\n$. The notation $\chi_\omega$ (respectively. $\chi_{\Omega \setminus \omega}$) denotes the characteristic function of $\omega$ (respectively. ${\Omega \setminus \omega}$). The second main difference arises from practical considerations: it is unrealistic from the point of view of applications to know the full graph of Dirichlet-to-Neumann. Therefore, we will assume that one has access to a single point in that graph. This non destructive testing problem is usually written from a numerical point of view as the minimization of a cost function: typically a least-square matching criterion. Many authors have investigated the steepest descent  method for this problem \cite{Kirsch,HettlichRundell1,ItoKunischLi,Pantz,AfraitesDambrineKateb} with the methods of shape optimization since the unknown parameter is a geometrical domain.

This work is devoted to the study of second order methods for this problem that has only be considered before for simplified models in \cite{EpplerHarbrecht,AfraitesDambrineEpplerKateb}. By introducing   second order methods, one aims to reach two distinct objectives.
\begin{itemize}
\item On  one hand, we  provide all the needed material to design
a Newton algorithm.  We will give differentiability results for
the state function and for the objective that  we have chosen to
study in this work. Nevertheless, we point out that the
discretization of a Newton method for this problem turns out to be
very delicate; this is why, in the present paper,  we will neither discuss  
about this problem nor present numerical examples. This topic
is actually the main objective of a work in progress. 
\item On the other hand, we analyze rigorously the well-posedness of the
optimization method. This is justified by the huge numerical
literature devoted to the numerical study of this question in the
field of inverse problems;  the numerical experiments insist on
the  ill-posedness of this problem. We will explain the
instability in the continuous settings in terms of shape
optimization. We show that the shape Hessian is not coercive -in
fact its Riesz operator is compact -- and this
explains the unstability of the minimization process.
\end{itemize}

Let us describe the precise problem under consideration and the notations. We consider a bounded domain $\Omega \subset \R^d$ ($d = 2$ or $3$) with a $\mathcal{C}^2$ boundary. It is filled with a material whose conductivity is $\sigma_{1}$ and with  an unknown inclusion $\omega$ in $\Omega$ of conductivity $\sigma_{2} \neq \sigma_{1}$. We search to reconstruct the shape of $\omega$ by measuring on $\partial \Omega$, the input voltage and the corresponding output current. In the sequel, we fix
$ d_0>0$ and consider inclusions $\omega$ such
that $\omega\subset\subset \Omega_{ d_0} =\{x \in \omega,~
d(x,\partial\Omega)> d_0\}$. We also assume that the boundary $\partial\omega$ is of class $\mathcal{C}^{4,\alpha}$.  The inverse problem arises when one has access to the normal vector derivative of the potential $u$ that solves \eqref{equation:etat}-\eqref{equation:potentiel:condition:limite:dirichlet}
when the conductivity distribution is defined by \eqref{definition:conductivite} . Assume that ones knows
\begin{equation}
\label{equation:potentiel:condition:limite:neumann}
\sigma_1\partial_\n u = g \text{ on } \partial \Omega,
\end{equation}
then the problem \eqref{equation:etat}-\eqref{equation:potentiel:condition:limite:dirichlet}-\eqref{equation:potentiel:condition:limite:neumann} is overdetermined. The electrical impedance tomography problem we  consider is to recover the shape of $\omega$ from the knowledge of the single Cauchy pair $(f,g)$.

In order to recover the shape of the inclusion $\omega$, an usual
strategy is to minimize a cost function. Many choices are
possible; however  it turns out that a Kohn and Vogelius type
objective leads to a minimization problem with nicer properties
than the least squares fitting approaches (we refer to
\cite{AfraitesDambrineKateb} for a comparison of different
objectives with order one methods and to
\cite{AfraitesDambrineEpplerKateb} for the case of a perfectly insulated inclusion). Therefore, we study such a
cost function  in this work.

Let us  define this criterion. Its distinctive feature is to
involve two state functions $\ud$ and $\un$: the state $\ud$
solves \eqref{equation:etat}-\eqref{equation:potentiel:condition:limite:dirichlet}
while $\un$ solves
\eqref{equation:etat}-\eqref{equation:potentiel:condition:limite:neumann}.
The  Kohn -Vogelius objective $\KV$ is then defined  as:
\begin{equation}
\label{definition:kohn:vogelius}
\KV(\omega)=\int_\Omega \sigma | \nabla(\ud-\un)|^2
\end{equation}
Let us sum up  the results of this paper concerning  the minimization of this objective.   We first prove differentiability results for the state $\ud$. In the sequel, we use the convention that a bold character denotes a vector. If $\h$ denotes a deformation field,  it can be written as $\h = \ct{h} + h_{n} \n$ on $\partial\omega$. 
Note also that in the following lines, $\n$ denotes the outer normal field to $\partial\omega$ pointing into $\Omega\setminus\overline{\omega}$. Hence, for $x\in \partial\omega$, we define, when the limit exists, $u^{\pm}(x)$ (resp. $(\partial_{n}u)^{\pm}(x)$) as the limit of $u(x \pm t\n(x))$ (resp. $\langle \nabla u(x\pm t \n(x),\n(x))$) when $t>0$ tends to $0$. Note that $\ct{h}$ is a vector while $h_n$ is a scalar quantity.

The admissible deformation fields have to preserve $\partial\Omega$ and the regularity of the boundaries: therefore the space of admissible fields is
\begin{equation*}
\mathcal{H} =\{\h \in \mathcal{C}^{4,\alpha}(\R^d, \R^d), Supp(\h) \subset  \Omega_{ d_0}\}.
\end{equation*}
The following result concerns the first order derivative of the state functions $\ud$ and $\un$. It was derived in \cite{HettlichRundell1, Pantz, AfraitesDambrineKateb}.
\begin{theorem}
\label{lem:derivee:etatD} Let $\Omega$ be an open smooth subset of $\R^d$
($d=2$ or $3$)  and let $\omega$
be an element of $\Omega_{d_0}$ with a boundary of class
$\mathcal{C}^{4,\alpha}$. Then the state functions $\ud$ and $\un$
are shape differentiable; furthermore  their shape derivative
$\ud'$ and $\un'$ belongs to
$\sH^1(\Omega\setminus\overline{\omega})\cup\sH^1(\omega)$ and
satisfy
\begin{equation}
\label{derivee:etatD}
\left\{
\begin{array}{rcl}
\Delta \ud'  & =  & 0 \text{ in }\Omega\setminus\overline{\omega} \text{ and  in } \omega,   \\
  \left[\ud'\right]& =  & h_n ~\cfrac{[\sigma]}{\sigma_{1}} ~ \partial_{\n}\ud^- ~\text{ on } \partial \omega, \\
  \left[\sigma \partial_{n} \ud'\right]&  = &  [\sigma] \divt{ h_n \nablat \ud}  \text{ on } \partial \omega,\\
  \ud' &=& 0 \text{ on }\partial\Omega.
\end{array}
\right.
\end{equation}
\begin{equation}
\label{derivee:etatN}
\left\{
\begin{array}{rcl}
\Delta \un'  & =  & 0 \text{ in }\Omega\setminus\overline{\omega} \text{ and  in } \omega,   \\
  \left[\un'\right]& =  & h_n ~\cfrac{[\sigma]}{\sigma_{1}} ~ \partial_{\n}\un^- ~\text{ on } \partial \omega, \\
  \left[\sigma \partial_{n} \un'\right]&  = &  [\sigma] \divt{ h_n \nablat \un}  \text{ on } \partial \omega,\\
  \partial\un' &=& 0 \text{ on }\partial\Omega.
\end{array}
\right.
\end{equation}
\end{theorem}
The main result of this work concerns the second order derivative.
It is given  is the following theorem.
\begin{theorem}
\label{theoreme:derivee:seconde:etatD} Let $\Omega$ be an open
smooth subset of $\R^d$ ($d=2$ or $3$) and let $\omega$ be an  element of $\Omega_{d_0}$ with a $\mathcal{C}^{4,\alpha}$ boundary.  Let ${\h}_1$ and ${\h}_2$ be
two deformation fields in $\mathcal{H}$.  Then the  state $\ud$
has a second order shape derivative $\ud'' \in
\sH^1(\Omega\setminus\overline{\omega})\cup\sH^1(\omega)$ that
solves
\begin{equation}
\label{derivee:seconde:etatD}
\left\{
\begin{array}{rcl}
\Delta \ud''  & =  & 0 \text{ in }\Omega\setminus\overline{\omega} \text{ and  in } \omega,   \\
  \left[\ud''\right]& =&  \left( h_{1,n} h_{2,n}H -\ct{{\h}_1}.(D\n\, \ct{{\h}_2}) \right) [\partial_{\n} \ud] - \left( h_{1,n}[\partial_{\n} (\ud)'_2] + h_{2,n}[\partial_{\n} (\ud)'_1] \right) \\
  && ~~~~~~+ \left(\ct{{\h}_1}.\nabla h_{2,n} +\ct{{\h}_2}.\nabla h_{1,n} \right) [\partial_{\n} \ud] \text{ on } \partial \omega, \\
  \left[\sigma \partial_{n} \ud'' \right]&  = & \divt{h_{2,n}\left[\sigma \nablat (\ud)'_1 \right] + h_{1,n}\left[\sigma \nablat (\ud)'_2 \right]  + \ct{h_1}. (D\n \,\ct{h_2})[\sigma \nablat \ud]}\\
 &&   ~~~~~~  - \divt{ (\ct{{\h}_1}.\nablat h_{2,n} +  \nablat h_{1,n}.\ct{h_2})  \left[\sigma \nablat \ud \right] }\\
 && ~~~~~~+ \divt{ h_{2,n} h_{1,n}  (2  D\n-H I) \left[\sigma \nablat  \ud \right]} \text{ on } \partial\omega,\\
   \ud'' &=& 0 \text{ on }\partial\Omega.
\end{array}
\right.
\end{equation}
\end{theorem}
Here, $(\ud)'_i$ denotes the first order derivative of $u$ in the direction of $h_i$ as given in \eqref{derivee:etatD}, $D\n$ stands for the second fundamental form of the manifold $\partial\omega$ and $H$ stands for the mean curvature of $\partial\omega$.
% In the bidimensional case, the jump of the flux can be written under the form
%\begin{eqnarray*}
% \left[\sigma \partial_{n} u''\right]& =& \divt{ \left(\ct{{\h}_1}.\ct{{\h}_2} - h_{1,n}.h_{2,n} \right)H [\sigma\partial_{\n} u]} + \divt{\left( h_{1,n} [\sigma \nablat u'_2] + h_{2,n} [\sigma \nablat u'_1] \right)}\\
 % && ~~~~~~ - \divt{\left(\ct{{\h}_1}.\nablat h_{2,n} +\ct{{\h}_2}.\nablat h_{1,n} \right) [\sigma \nablat  u]  } \text{ on } \partial \omega.
%\end{eqnarray*}
%\end{theorem}
The twin result concerning $\un$ is an easy adaption of Theorem \ref{theoreme:derivee:seconde:etatD}.
Once the differentiability of the state function has been established, one can consider the objectives.
In \cite{AfraitesDambrineKateb}, we have shown the first order result.

\begin{theorem}
\label{derivee:premiere:kohn:vogelius}
Let $\Omega$ be an open smooth subset of $\R^d$ ($d=2$ or $3$) and let $\omega$ be an element of $\Omega_{d_0}$ with a $\mathcal{C}^{4,\alpha}$ boundary.  Let ${\h}_1$ and ${\h}_2$ be two deformation fields in $\mathcal{H}$. The  Kohn-Vogelius objective is differentiable with respect to the shape and its derivative in the direction of a deformation field $\h$ is given by:
\begin{equation}
\label{derivee:premiere}
D\KV(\omega)\h=  [\sigma] \int_{\partial\omega} \left[ \cfrac{\sigma_1}{\sigma_2} \left( |\partial_\n \ud^+|^2 - |\partial_\n \un^+|^2 \right) + |\nablat \ud|^2 -|\nablat \un|^2 \right] h_n.
\end{equation}
\end{theorem}

We now give  the second-order derivative  of the Kohn and Vogelius criterion. 

\begin{theorem}
\label{derivee:seconde:kohn:vogelius} Let $\Omega$ be an open smooth
subset of $\R^d$ ($d=2$ or $3$) and $\omega$ be an element of $\Omega_{d_0}$ with a $\mathcal{C}^{4,\alpha}$ boundary.  Let ${\h}_1$ and ${\h}_2$ be
two deformation fields in $\mathcal{H}$. The  Kohn-Vogelius
objective is twice differentiable with respect to the shape and
its second derivative in the directions $\h_1$ and $\h_2$ is given
by:
\begin{equation}
\label{second:derivative}
\begin{split}
D^2\KV(\omega)(\h_1,\h_2)&=\int_{\partial \omega} \left[ \sigma | \nabla v |^2 \right]\left(\ct{\h_1}.\nablat(h_{2,n})+ \ct{\h_2}.\nablat(h_{1,n})- \ct{h_{2}} .(D\n\, \ct{\h_1}) \right)\\
&-\int_{\partial \omega}\partial_\n \left(\left[\sigma | \nabla  v |^2\right] \right) h_{1,n} h_{2,n}+ 2 \big[\sigma \nabla v .( h_{1,n}\nabla v'_2 +h_{2,n} \nabla v'_1 )\big]\\
&+\int_{\partial \omega} \left[\sigma \left(\partial_\n  (\un)'_1 v'_2 +  \partial_\n  (\un)'_2 v'_1 -  \partial_\n v'_1 (\ud)'_2-   \partial_\n v'_2  (\ud)'_1 \right)\right]\\
&+2  \int_{\partial \omega}  v \left[ \sigma \partial_\n (\un)''_{1,2} \right] - \sigma_1 \partial_\n v^+ \left[(\ud)''_{1,2}\right]
\end{split}
\end{equation}
where we have set $v=\ud-\un$.
\end{theorem}
To investigate the properties of stability of this cost function,
we are led to  consider an admissible inclusion $\omega^*$ to
solve both
\eqref{equation:etat}-\eqref{equation:potentiel:condition:limite:dirichlet}
and
\eqref{equation:etat}-\eqref{equation:potentiel:condition:limite:neumann}
in order to obtain the corresponding measurements $f^*$ and $g^*$.
It is obvious that the domain $\omega^*$ realizes
the absolute minimum of the criterion $\KV$ since,  by
construction, we can write  $\ud=\un $ in $\Omega$ and hence
$\KV(\omega^*) =0$. We will check that the Euler equation
\begin{equation*}
D\KV(\omega^*)(\h)=0,
\end{equation*}
holds.  We will also prove that
\begin{equation}
\label{positivite:hessienne:forme:minimum:absolu}
D^2\KV(\omega^*)(\h,\h)=\int_{\Omega}\sigma |\nabla v'|^2 .
\end{equation}
Moreover, if $h_n \ne 0$, then 
$D^2\KV(\omega^*)(\h,\h)>0$ holds. Nevertheless,
\eqref{positivite:hessienne:forme:minimum:absolu} does not mean
that the minimization problem is well-posed. In fact, it is the
following theorem that explains the instability of standard
minimization algorithms.
\begin{theorem}
\label{compactness:D2KV} Assume that $\omega^*$ is a critical
shape of $\KV$ for which the additional condition $\un = \ud$
holds. Then the Riesz operator corresponding to $D^2\KV(\omega^*)$
defined from $ \sH^{1/2}(\partial\omega^*)$ with values in
$\sH^{-1/2}(\partial\omega^*)$ is compact. Moreover, the
minimization problem is severely ill-posed in the following sense:
if the target domain is $\mathcal{C}^\infty$ and if $\lambda_n$
denotes the $n^{th}$ eigenvalue of $D^2\KV(\omega^*)$, then
$\lambda_n = o(n^{-s})$ for all $s>0$.
\end{theorem}

Theorem \ref{compactness:D2KV} has two main consequences. First,
the shape Hessian at the global minimizer is not coercive. This
means that this minimizer may not be  a local strict minimum of the
criterion. Moreover, the criterion provides no control of the
distance between the parameter $\omega$ and the target $\omega^*$.
The second consequence  concerns any numerical scheme used to
obtain this optimal domain $\omega^*$. One has to face this
difficulty and this  explains why frozen Newton or Levenberg-Marquard schemes have been  used  to  solve numerically this problem \cite{HettlichRundell1,AfraitesDambrineKateb}.

The paper is organized as follows. In a first section, we state some
preliminary results. Some are well known facts in shape
optimization and will be recalled without proof for the sake of
readability. Some of them (e.g the derivatives of a
Laplace-Beltrami operator and the tangential regularity of the
solution to
\eqref{equation:etat}-\eqref{equation:potentiel:condition:limite:dirichlet}
along the discontinuity of the  conductivity distribution) are
less known and will be proved thanks to potential layer methods.
Hence we will tackle the computations  in Section
\ref{section:preuve:differentibilite} that we consider as  the
core of this work : it is essentially devoted to prove Theorem
\ref{theoreme:derivee:seconde:etatD}. After a first part  where we
prove the existence of a second order derivative for the state, we
propose two distinct methods to find the boundary value problem
solved by this second order derivative. The first method
(subsection \ref{section:preuve:differentibilite:etape2}) follows
the lines of classical proofs of shape differentiability by
differentiating the weak formulation of problem
\eqref{equation:etat}-\eqref{equation:potentiel:condition:limite:dirichlet}
and interpreting the result in terms of differential operator and
boundary conditions. The alternative method (subsection
\ref{section:preuve:differentibilite:etape3}) consists in a direct
differentiation of the boundary conditions. Finally, Section
\ref{section:derivee:seconde:critere} is devoted to the analysis
of the criterion, we establish Theorem
\ref{derivee:seconde:kohn:vogelius} and Theorem
\ref{compactness:D2KV}. We will present their consequences on  the
stability of critical shapes.

%%%%%%%%%%%%%%%%%%%%%%%%%%%%%
\section{Preliminary results. }
\label{annexe}

\subsection{Elements of shape calculus}
% \subsection{Notations and geometrical facts.}

Before entering the proof of Theorem \ref{theoreme:derivee:seconde:etatD}, we  recall without proof
some basic facts from shape optimization (see \cite{HenrotPierre} for references).
Let $\h$ be a deformation field in $\mathcal{C}^2(\Omega,\R^d)$ with $\|\h\|_{\mathcal{C}^2} <1$. We set $T_t(h,.) =Id +t \h$ and denote by $\Omega_t$ the
transported domain $\Omega_t=T_t(\Omega)$. To avoid heavy notations, we will misuse the notation $T_t$ instead of $T_t(h,.).$

\textbf{Material and shape derivatives.}  Classically, in mechanics of continuous media, the material derivative is defined  as being a positive limit. In our context,
 for any vector field $\h \in \mathcal{H}$, we define the material derivative of the domain functional
$y=y(\Omega)$ at $\Omega$ in an admissible  direction  $\h$
 as the limit
 \begin{equation}
\dot{y}(\Omega;\h)=\lim_{t\rightarrow 0}\cfrac{y(\Omega_t)\circ T_t-y(\Omega)}{t},
 \end{equation}

\par
\noindent
Similarly, one can define the material derivative
$\dot{y}(\partial\Omega,\h)$ for any domain functional $y=y(\partial\Omega)$ which
depends on $\partial\Omega$. Another kind of derivative occurs : it is
called the shape derivative of $y(\Omega,\h)$. It is viewed as a
first local variation. Its definition is given by the following
\begin{definition}
The shape derivative $y'=y'(\Omega;\h)$ of a functional $y(\Omega)$
at $\Omega$ in the direction of a vector field $\h$ is given  by
\begin{equation}
y'=\dot{y}- \h.\nabla y.
\end{equation}
\end{definition}
For more details on these derivations, the reader can consult
\cite{SokolowskiZolesio,HenrotPierre}.

\textbf{Elements of tangential derivatives.}
We will need in the sequel to manipulate the tangential
differential operators on a manifold. For the reader's
convenience, we recall from \cite{DelfourZolesio,HenrotPierre}
some  definitions and also some useful rules of calculus.
\begin{definition}
The tangential divergence  of a vector field $\mathbf{V}\in
C^1(\R^d,\R^d)$  is given by
\begin{equation}
\label{definition:divergence:tangentielle}
\divt{\mathbf{V}}=\div{\mathbf{V}} -  D\mathbf{V}.\n.\n,
\end{equation}
where the notation $D\mathbf{V}$ denotes  the Jacobian matrix of
$\mathbf{V}$. When the vector $\mathbf{V} \in
C^1(\partial\Omega,\R^d)$ is defined on $\partial\Omega$, then the
following notation is used to define the tangential divergence
\begin{equation}
\divt{\mathbf{V}}=\div{ \tilde{\mathbf{V}}}- (D\tilde{\mathbf{V}}.\n).\n,
\end{equation}
where $\tilde{\mathbf{V}}$ stands for an arbitrary $C^1$ extension of $\mathbf{V}$ on an open neighborhood of $\partial\Omega$.
\end{definition}

We introduce now, the notion of  tangential gradient $\nablat$ of
any smooth scalar function $f$ in $\mathcal{C}^1(\partial
\Omega,\R^d)$.
\begin{definition}
Let an element $f\in \mathcal{C}^1(\partial \Omega,\R^d)$ be given
and let $\tilde f$ be an extension of $f$ in the sense that
$\tilde f \in \mathcal{C}^1(U)$ and $\tilde f|_{\partial \Omega}=f$ and
where $U$ is an open neighborhood of $\partial \Omega$. Then the
following notation is used to defined the tangential gradient
\begin{equation}
 \nablat f= \nabla \tilde f|_{\partial \Omega} -  \nabla \tilde f.\n ~\n \text{ on }\partial\Omega.
\end{equation}
\end{definition}
The details for the existence of such an extension can be found in
\cite{DelfourZolesio}. Let us remark that these definitions do not
depend on the choice of the  extension.
Furthermore, one can show the important relation
\begin{equation}
\int_{\partial\Omega } \nablat f .\mathbf{F} ~= - \int_{\partial\Omega} f~\divt{ \mathbf{F}},
\end{equation}
for all elements  $f \in \mathcal{C}^1(\partial\Omega)$ and all vector
fields $F \in \mathcal{C}^1(\partial\Omega,\R^d)$ satisfying
$F_n=\langle F,n \rangle=0$.

\textbf{Integration by parts on $\partial \Omega$.}
In general, the condition above $F_n=0$ is not always satisfied.
We are then led to find another formula to extend the formula in
the general case. The extension of this integration by parts
formula to fields with a normal vector component involves
curvature.

First, we point out that the curvature is connected to the normal
vector via the tangential divergence operator. Recall that the
mean curvature of $\partial\Omega$ is defined as $
H=\text{div}_\tau(\n)$. Making use of the form of
$\text{div}_\tau(\n)$ on the boundary, one shows straightforwardly the following statement.

\begin{proposition}
Let $\Omega$ be an
open subset of $\R^3$ with a $\mathcal{C}^2$ boundary. For any
unitary extension ${\cal N}$ of $n$ on a neighborhood of $\partial
\Omega$, one has
\begin{equation*}
\div{ \mathcal{N}}=H\text{ on } \partial\Omega.
\end{equation*}
Assume that  the manifold $\partial\Omega$
has no borders. If $\mathbf{F}\in
\sH^2(\partial\Omega)^3$ and $f\in \sH^2(\partial\Omega)$, then we
have
\begin{equation}
\label{formule:integration:par:parties}
\int_{\partial\Omega} \nabla f .\mathbf{F} +f \divt{\mathbf{F}} = \int_{\partial\Omega}\left(
{\nabla f . \n}+H f \right) \mathbf{F}.\n.
\end{equation}
\end{proposition}
%\subsection{Second order derivative}
We assume now that the domain $\Omega$ has a $\mathcal{C}^3$
boundary. The simplest second-order derivative is the Laplace
Beltrami operator; it is defined as follows (see
\cite{SokolowskiZolesio,DelfourZolesio,HenrotPierre}) thanks to
the following usual chain rule.
\begin{definition}
Let $f \in H^2(\partial \Omega)$. The Laplace-Beltrami $\Deltat$
of $f$ is defined as follows
\begin{equation}
\Deltat{f}=\divt{\nablat f}.
\end{equation}
\end{definition}
There is a  relation connecting the Laplace operator and the
Laplace-Beltrami operator. Let us denote by $\partial^2_{nn} f=  (D^2f.\n).\n$ where $D^2f$ stands for the Hessian of $f$.
\begin{proposition}
Let $\Omega$ be a domain with  a boundary $\partial\Omega$ of class $\mathcal{C}^3$.
For all functions $f\in H^3(\Omega)$, it holds
\begin{equation}
\label{lien:laplacien:laplace-beltrami} \Delta f=\Deltat f + H
\partial_\n f+\partial^2_{nn} f, \text{ on } \partial\Omega.
\end{equation}
\end{proposition}

% \subsection{Some useful derivatives. }

We need to compute  shape and material derivative of special
vector fields: the outer unit normal vector $\n$,  the tangential
gradient and  the Laplace-Beltrami operator applied to a function. While  the derivative of
the normal vector  is obtained by a straightforward calculus,   we
have to  transport from $\partial\Omega_t$ to $\partial\Omega$ the
Laplace-Beltrami operator and the tangential gradient in order to
compute the other  derivatives.

\emph{Derivatives of the normal vector.}
We  describe  the material and shape  derivatives of the normal vector.
We will denote by $\n$  the gradient of the signed distance to  $\partial\Omega$.
This is  an unitary extension of the
 unitary normal vector $\n$ at $\partial\Omega$ which is smooth in the vicinity of $\partial\Omega$.
 This extension furnishes a symmetric Jacobian $D\n$ that satisfies $D\n\,\n=0$ on $\partial\Omega$.
 The direction $\h$ will be supposed to be in
 $\mathcal{C}^2(\R^d,\R^d)$ or in $\mathcal{C}^{2}(\partial\Omega,\R^d)$.

\begin{proposition}\label{derivative:normal}
The material derivative $\dot{\n}$ of the normal vector $\n$ at $\Omega$ in
the direction of a vector field $\h\in
\mathcal{C}^1(\R^d,\R^d)$ is given by
\begin{equation*}
\dot{\n}=-\nablat(\h.\n)+D\n \,\ct{\h},
\end{equation*}
where $\ct{\h}=\h-\h. \n~\n$.
\end{proposition}
Concerning its shape derivative defined as $\n' = (\partial_t \n_t) |_{t=0}$ where $\n_t$ is any smooth unitary extension of $\n$  to  $\partial\Omega_t$,  we obtain.\begin{proposition}
The shape boundary $\n'$ in the direction of $\h$ is given by
\begin{equation*}
\n'=-\nablat(\h.\n).
\end{equation*}
\end{proposition}

\emph{Derivative of the tangential gradient.} For  $f\in \sH^{3}(\partial\Omega)$, we  compute the material derivative of $\nablat f$. %and $\Deltat f$ . Concerning $\nablat f$.
We  first compute the difference $\dot{\overline{\nablat f}}-\nabla \dot{f}$.

\begin{proposition}
\label{derivee:materielle:gradient:tangentiel} For all functions
$f\in \mathcal{C}^2(\mathbb{R}^3)$ and directions $\h \in
C^2(\partial\Omega,\mathbb{R}^3)$, one has
\begin{equation*}
\dot{\overline{\nablat f}}=\nabla \dot{f}+(D^2 f \h)_{\tau} - \nabla f.\n~\dot{n}-\nabla f.\dot{\n}~\n
\end{equation*}
\end{proposition}

\begin{proofof}{Proposition \ref{derivee:materielle:gradient:tangentiel}.} We differentiate  $\nabla f$ and $\nabla f.\n~\n$ and obtain
\begin{equation*}
\dot{\overline{\nabla f} }= \nabla f' + D^2 f \h
\end{equation*}
while
\begin{equation*}
\dot{\overline{\nabla f.\n~\n}}=\nabla f.\dot{\n}~\n + \nabla f.\n ~ \dot{\n}+ \nabla f'.\n ~ \n +(D^2 f \h).\n~\n.
\end{equation*}
The two former equations  give the desired result.
\end{proofof}

\emph{Derivative of the Laplace-Beltrami operator.}
Now, we want to compute the material derivative $\dot{\overline{\Deltat f}}$. We
begin  to study how to transport the Laplace-Beltrami operator
when one works on  $\partial\Omega_t$.  Let $\Deltatt$ denote the Laplace-Beltrami operator on the manifold $\partial\Omega_t$. To compute the derivative of a Laplace-Beltrami operator, we need the following proposition that we quote  from \cite{SokolowskiZolesio}.
\begin{proposition}
\label{prop:soko:zol} Let $f \in \sH^{5/2}(\R^d)$, then
\begin{equation}
\int_{\partial\Omega}  \left[ \left(\Deltatt f\right) \circ
T_t~\gammat(t) \right] \phi~ = - \int_{\partial\Omega} \left[
C(t)~\Big( \nabla(f\circ T_t)- (B(t)~\n).\nabla(f\circ
T_t)\Big)\right] .\nabla \phi,\,\,\forall \phi \in \mathcal{D}(\R^d).
\end{equation}
\end{proposition}
In the former proposition, we set
\begin{equation}
\begin{array}{lll}
\gamma(t)&=&\det{D T_t}, \\
\gammat(t)&=&\gamma(t) \| (DT_t^{-1})^T.\n\|_{\R^d},\\
B(t)&= &\cfrac{D(T_t^{-1})(D(T_t^{-1})^T} {\| (DT_t^{-1})^T.\n\|_{\R^d}^2},\\
C(t)&=&\gammat(t) D(T_t^{-1})(D(T_t)^{-1})^T.
\end{array}.
\end{equation}
A straightforward computation gives
\begin{equation}
\begin{array}{lll}
\gamma'(0)&=&\divt{\h},\\
\gammat'(0)&=&\divt{\h} = \divt{\ct{\h}}+H \h_n ,\\
B'(0)&= &2( D{\h} \n).\n I-( D{\h}+( D{\h})^T),\\
C'(0)&=& \divt{\h} I - ( D{\h}+( D{\h})^T).
\end{array}.
\end{equation}

\begin{theorem}
\label{derivee:materielle:laplace:beltrami}
Let $f \in {\cal{D}}(\R^d)$. The material derivative of $\Deltat f$ in the
direction $\h$ is given by
\begin{equation}
\label{formule:derivee:materielle:laplace:beltrami}
\begin{split}
\dot{\overline{\Deltat f }} = \Deltat \dot{f} + &\nablat f.\nablat\left[\divt{\ct{\h}}\right] + \nablat(H \h_n).\nablat f - \divt{\left(\left( D{\h}+( D{\h})^T \right)\nablat f \right)_{\tau}}
\end{split}
\end{equation}
\end{theorem}

\begin{proofof}{Theorem \ref{derivee:materielle:laplace:beltrami}}:
Formula \eqref{formule:derivee:materielle:laplace:beltrami} is shown in a weak sense. For each test function $\phi \in \mathcal{C}^\infty(\partial\Omega)$, there exists an extension $\phiet \in \mathcal{D}(\R^d)$ such that $\partial_\n \phiet =0$; this can be done by extending $\phi$ as a constant along the orbits of the gradient of the signed distance function to $\partial\Omega$ and the use of a cut-off function. For $f \in \mathcal{D}(\R^d)$, we set
\begin{equation*}
A(t) = \int_{\partial\Omega} \cfrac{(\Deltatt f)\circ T_t-\Deltat f }{t}~\gammat(t) ~\phi.
\end{equation*}
After an integration by parts on $\partial\Omega$, we obtain:
\begin{eqnarray*}
A(t)& = &\int_{\partial\Omega}  \cfrac{1-\gammat(t)}{t} ~(\Deltatt f)\circ T_t ~ \phi
+ \int_{\partial\Omega} \frac{\gammat}{t} \left( (\Deltatt f)\circ T_t  \phi + \frac{1}{t}\nablat f. \nablat \phi \right),\\
 & = & \int_{\partial\Omega}  \cfrac{1-\gammat(t)}{t} ~(\Deltatt f)\circ T_t ~ \phi\\
 &&~~~
 +\int_{\partial\Omega} \frac{1}{t} \left( \left[ \nablat f -C(t) \nabla \left( f\circ T_t \right) \right].\nabla \phiet+ \left[(B(t) \n. \nabla(f\circ T_t)\right] C(t)~\n.\nabla \phiet \right).
\end{eqnarray*}
Since $\partial_\n \phiet=0$ and $C(0) = I$, we get
\begin{equation*}
A(t) =  \int_{\partial\Omega}  \cfrac{1-\gammat(t)}{t} ~(\Deltatt f)\circ T_t ~ \phi
+ \int_{\partial\Omega}    \frac{ \nablat (f - f\circ T_t) }{t}  .\nablat \phiet  + \int_{\partial\Omega}    \frac{ C(0)-  C(t) }{t} \nabla (f\circ T_t) .\nablat \phiet.
\end{equation*}
When $t \rightarrow 0$, it then comes
\begin{eqnarray*}
\int_{\partial\Omega} \dot{\overline{\Deltat f}} \phi&=& -\int_{\partial \Omega} \gammat'(t) \Deltat f \phi + \nablat \dot{f}.\nablat \phi + \left(C'(0).\nabla f \right).\nablat \phi,\\
&=& \int_{\partial\Omega} \left( \Deltat \dot{f} - \divt{\h} \Deltat f \right) \phi + \left(  D{\h}+( D{\h})^T - \divt{\h} I\right)\nabla f .\nablat \phi,\\
&=& \int_{\partial\Omega} \left[ \Deltat \dot{f} - \divt{\h} \Deltat f  + \divt{\divt{h} \nablat f} - \divt{\left( \left(  D{\h}+( D{\h})^T\right) \nabla f \right)_{\tau}} \right]  \phi.
\end{eqnarray*}
Expanding the double divergence term, we  obtain:
\begin{equation*}
\dot{\overline{\Deltat f}} = \Deltat \dot{f}  + \nablat f.\nablat \divt{\h} - \divt{\left(\left( D{\h}+( D{\h})^T \right)\nabla f \right)_{\tau}}.
\end{equation*}
In order to explicit these derivatives, we let appear the curvatures of $\partial\Omega$ by means of
\begin{equation*}
 \nablat f.\nablat \divt{\h} =  \nablat f.\nablat \left[ \divt{\ct{\h}} + H \h_n\right],
\end{equation*}
%and also by splitting $\nabla f$ into its tangential  and normal vector componants
%\begin{equation*}
%\divt{\left(\left( D{\h}+( D{\h})^T \right)\nabla f \right)_{\tau}} = \divt{\left(\left( D{\h}+( D{\h})^T \right)\nablat f \right)_{\tau}} + \divt{\partial_\n f\left(\left( D{\h}+( D{\h})^T \right)\n \right)_{\tau}}.
%\end{equation*}
%where thanks to the symetry of $D\n$, a duality argument leads to
%\begin{equation*}
% \divt{\left(\left( D{\h}+( D{\h})^T \right)\nablat f \right)_{\tau}} = 2 \divt{ D\ct{\h}\nablat f + \h_n D\n \,\nablat f }.
%\end{equation*}
and this ends the proof of the theorem \eqref{formule:derivee:materielle:laplace:beltrami}.
\end{proofof}
%
%\begin{remark}
%In the forthcoming computations, the functions $f$ to which this theorem will be applied have a vanishing flux that is to say $\partial_\n f=0$. In that particular case \eqref{formule:derivee:materielle:laplace:beltrami} simplifies in
%\begin{equation*}
%\dot{\overline{\Deltat f }} = \Deltat \dot{f} + \nablat \divt{\ct{\h}}.\nablat f + \nablat(H \h_n).\nablat f - \divt{\left(\left( D{\h}+( D{\h})^T \right)\nablat f \right)_{\tau}}
%\end{equation*}
%\end{remark}

%\begin{remark}
%If the deformation field is tangential $\h=\ct{\h}$, then \eqref{formule:derivee:materielle:laplace:beltrami}
%simplifies and one has
%\begin{equation}
%\dot{\overline{\Deltat f }} = \Deltat \dot{f} + \nablat \divt{\ct{\h}}  -2 \divt{D\ct{\h}  \nablat f }  + \divt{ \partial_\n f  \left( \left( D\ct{\h} +(D\ct{\h})^T \right) \n \right)_{\tau}}.
%\end{equation}
%\end{remark}

% For all $s\in(1/2,4]$, the operator  $K_\omega: \sH^{s}(\partial \omega) \rightarrow \sH^{s}(\partial \omega)$ is compact. Hence, the Fredholm alternative provides the existence and uniqueness of a solution for all $s\in(1/2,4)$ in $\sH^s(\partial\omega)$ for each $s\in(1/2,4]$. This common solution is $(\ud)^+_{|\partial\omega}$.
%\end{proofof}
%Corollary \ref{regularite:tangentielle} provides regularity for the solutions of  \eqref{systemederivee:ud}; this is needed  to justify the computations of the shape derivatives.
%%%%%%%%%%%%%%%%%%%%%%%%%%%%%
\section{Existence of the second order derivative of the state.
Proof of Theorem \ref{theoreme:derivee:seconde:etatD}.}
%%%%%%%%%%%%%%%%%%%%%%%%%%%%%
\label{section:preuve:differentibilite}

The section is devoted to prove Theorem
\ref{theoreme:derivee:seconde:etatD}. We follow the  usual
strategy to derive existence in shape optimization. In section
\ref{section:preuve:differentibilite:etape1}, we will write the
weak formulation of the problem, then transport it on the
reference domain, pass to the limit and obtain existence of the
material derivative. In a second time, we will seek a boundary
value problem solved by the material derivative. This will provide
a characterization of the second order shape derivative. Two
strategies, that we will detail,  are possible: the first one
explored in section \ref{section:preuve:differentibilite:etape2}
consists in working on the variational formulation while the
second one uses the tangential differential calculus by
differentiating the boundary conditions. This last approach  will
be  presented in section \ref{section:preuve:differentibilite:etape3}.  The computations that will be made in subsections \ref{section:preuve:differentibilite:etape2} and \ref{section:preuve:differentibilite:etape3} require some regularity of the traces of the state $\ud$ on the interface of discontinuity $\partial\omega$. For the sake of readability, we postponed in subsection \ref{section:preuve:differentibilite:etape4} all the needed justifications.

\subsection{Preliminary  results.}

In the sequel, we will use some technical formulae. To preserve the readability of the proof of the main result, we state them in this paragraph. The tools needed  for proving these results can be found in \cite{SokolowskiZolesio}. Given a smooth  vector field $\h$, we denote
\begin{equation*}
A_{\h}=D\h+D\h^T-\div{\h}I
\end{equation*}
We begin with the following formula.
\begin{lemma}
\label{fait:technique0}
It holds:
\begin{equation}\label{formule:magique1}
\nabla u . A_{\h} \nabla v = \nabla (\h.\nabla u ). \nabla v + \nabla (\h.\nabla v )\nabla u -
\div{(\nabla u . \nabla v ) \h}.
\end{equation}
\end{lemma}
Given two smooth vector fields $\h_1$ and $\h_2$, we set
\begin{equation}
\label{definition:AA}
\AA = D{\h}_2 A_{{\h}_1}+A_{{\h}_1} D{{\h}_2}^T-A_{{\h}_1}\div{{\h}_2} -(A_{{\h}_1})'({\h}_2),
\end{equation}
and
\begin{equation*}
b=({\h}_2.\nabla{u})A_{{\h}_1}\nabla{v}+({\h}_2.\nabla{v})A_{{\h}_1}\nabla{u}-((A_{{\h}_1}\nabla{u}).\nabla{v}){\h}_2.
\end{equation*}
Here, the notation $(A_{{\h}_1})'({\h}_2)$ stands  for  the matrix defined by its elements  $$((A_{{\h}_1})'({\h}_2))_{k,l}=\nabla(((A_{{\h}_1})')_{k,l}).\h_2$$
\begin{lemma}
\label{fait:technique1}
One has:
\begin{equation}
\label{formule:magique:2}
\nabla{u}.\AA\nabla{v}=\div{b}-({\h}_2.\nabla{u})\div{(A_{{\h}_1}\nabla{v})}-({\h}_2.\nabla{v})\div{(A_{{\h}_1}\nabla{u})}.
\end{equation}
\end{lemma}

%For the proof of Lemma \ref{fait:technique0} and Lemma \ref{fait:technique1}, we refer to .
We  need the following crucial result
\begin{lemma}
\label{fait:technique2}
If $u$ is harmonic  then
%$$ \div{A_{{\h}_1}\nabla{u}-\nabla{({\h}_1.\nabla{u})}}=0 $$

\begin{equation}\label{formule:harmonique}
\div{A_{{\h}_1}\nabla{u}}=\Delta({\h}_1.\nabla u).
\end{equation}
\end{lemma}

\begin{proofof}{Lemma \ref{fait:technique2}}
For any harmonic function $u$  in $\Omega$ and for every test function $\phi \in \mathcal{D}(\Omega)$, we can write
\begin{equation*}
\int_{\Omega} \nabla \dot{u} \nabla \phi= \int_{\Omega}  A_{\h} \nabla u  \nabla \phi
\end{equation*}
then
\begin{equation*}
\int_{\Omega} \Delta \dot{u}~~  \phi= \int_{\Omega} \div{ A_{\h} \nabla u }~~ \phi
\end{equation*}
Since $\dot{u}=u'+\h.\nabla u$ and since $u'$ is harmonic in  $\Omega$,  we obtain the result.
\end{proofof}
%%%%%%%%%%%%%%%%%%%%%%%%%%%%%
\subsection{Proof of existence of the second order derivative.}
\label{section:preuve:differentibilite:etape1}
%%%%%%%%%%%%%%%%%%%%%%%%%%%%%
We follow Hettlich and Rundell \cite{HettlichRundell2} and Simon \cite{Simon} to  define the second order derivative of an operator with respect to a
domain. We compute the second derivative by considering two  admissible deformations $\h_1,\h_2 \in \mathcal{H}$
that will describe the small variations of $\partial\omega$. Simon shows that the second derivative $F''(\partial \omega;\h_1,\h_2)$ of $F(\partial \omega)$  is defined as a bounded bilinear operator satisfying
\begin{equation}
F''(\partial \omega;\h_1,\h_2)=\left(F'(\partial \omega;\h_1)\right)'\h_2-F'(\partial \omega;D\h_1\, \h_2)
\end{equation}
For more details, the reader can consult the appendix in page 613 of \cite{HettlichRundell2}.

Let us begin the proof. Let ${{\h}_1},{{\h}_2} \in \mathcal{H}$ be two vector fields. The direction
${\h}_1$ being fixed, we consider $\dot{u}_{1,{\h}_2}$ the
variation of ${\dot{u_1}}$ with respect to the direction ${\h}_2$.
We recall from \cite{AfraitesDambrineKateb} that the material
derivative ${\dot{u_1}}$ of $u$ in the direction ${\h}_1$ satisfies
\begin{equation*}
 \label{limite:difference:variationelle}
\forall v \in H^1_0(\Omega),~~\int_{\Omega} \sigma \nabla \dot{u}_1.\nabla v=   \int_{\Omega}\sigma  \nabla u.{A}_{{\h}_1}\nabla v.
\end{equation*}

 Let $\phi_2 :\Omega
\mapsto\Omega$ be  the diffeomorphism defined by
$\phi_2(x)=x+{\h}_2(x)$ and we set $\psi_2=\phi_2^{-1}$. Setting $\omega_{{\h}_2}=\left\{x+{\h}_2(x),~x \in \omega\right\}$, $\Omega_{{\h}_2}=\left\{x+{\h}_2(x),~x \in \Omega\right\}=\Omega$ and $\sigma_{\h_2}=\sigma\circ \phi_2$, we get
\begin{equation}
\label{equation:domaine:fixe}
\int_{\Omega_{{\h}_2}}\sigma_{{\h}_2}
\nabla\dot{u}_{1,{\h}_2}.\nabla{v} =\int_{\Omega_{{\h}_2}}\sigma_{{\h}_2}
\nabla u_{{\h}_2}.A_{{\h}_1}\nabla{v}
\end{equation}
where $u_{{\h}_2}$ is the solution of the original problem with
$\omega_{\h_2}$ instead of $\omega$. Making the change of
variables $x=\phi_2(X)$, we get the integral identity on the fixed
domain $\Omega$ :
\begin{equation}\label{changing_variable}
\int_{\Omega}\sigma
\nabla\widetilde{\dot{u}}_{1,{\h}_2}.\left(D\psi_2 (D\psi_2)^T det(D\phi_2)\right)\nabla{v}
=\int_{\Omega}\sigma \nabla
\tilde{u}_{{\h}_2}.\left(D\psi_2\widetilde{A_{{\h}_1}}
 (D\psi_2)^T~det(D\phi_2)\right)\nabla{v}
\end{equation}
with the notations $\tilde{u}=u\circ\phi_2$ and
$\widetilde{A_{{\h}_1}}=A_{{\h}_1}\circ\phi_2$.
Since the material derivative
$\dot{u}_1$ of $u$ with respect to the direction  ${\h}_1$ satisfies
\begin{equation*}\label{material}
\int_{\Omega}\sigma
\nabla{\dot{u}_1}.\nabla{v} =\int_{\Omega}\sigma
\nabla{u}.A_{{\h}_1}\nabla{v},
\end{equation*}
the difference of \eqref{equation:domaine:fixe} and \eqref{changing_variable}  gives
\begin{equation*}
\begin{split}
\int_{\Omega}\sigma
\nabla{\left(\widetilde{\dot{u}}_{1,{\h}_2} - \dot{u}_1\right)}&.\nabla{v} =\int_{\Omega}\sigma\nabla{\widetilde{\dot{u}}}_{1,{\h}_2}
.\left(I-D\psi_2 (D\psi_2)^T~det(D\phi_2)\right)\nabla{v} \\
&+ \int_{\Omega}\sigma\nabla{\tilde{{u}}}_{{\h}_2} .\left(D\psi_2\widetilde{A_{\h_1}}
(D\psi_2)^T~det(D\phi_2)-A_{{\h}_1}\right)\nabla{v} +\int_{\Omega}(\nabla{\tilde{u}_{{\h}_2}}-\nabla{u}).A_{{\h}_1}\nabla{v}.
\end{split}
\end{equation*}

We quote from \cite{Kirsch} and \cite{HettlichRundell2} the following   asymptotic  formulae
\begin{eqnarray*}
\|- \div{\h_i} \|_{\infty}  &=&  O(\|\h_i \|_{\mathcal{C}^2}^2) ,\\
\|  D\psi_i (D\psi_i)^T  \det(D\phi_i)  - I+A_{\h_i} \|_{\infty}  &=&  O(\| \h_i \|_{\mathcal{C}^2}^2),\\
\|   D\psi_2 \widetilde{A_{{\h}_1}} (D\psi_2)^T \det(D\phi_2) -A_{{\h}_1} +  D{\h_2} A_{{\h}_1} + A_{{\h}_1} ( D{\h_2})^T- \div{{{\h}_2}} A_{{\h}_1}-(A_{{\h}_1})'({{\h}_2}) \|_{\infty}& =&  O(\| {{\h}_2} \|_{\mathcal{C}^2}^2).
\end{eqnarray*}
 Making the adequate substitutions, we easily check that the material
derivative of $\dot{u}_{1}$ with respect to ${\h}_2$ exists. This
derivative, denoted  by $\ddot{u}_1$, satisfies
\begin{equation}
\label{derivee:seconde:materielle:faible}
\int_\Omega \sigma \nabla{\ddot{u}_{1}}.\nabla{v}~dx=  \int_\Omega \sigma \left[ \nabla{\dot{u}_1}.A_{{\h}_2}\nabla{v} + \nabla{\dot{u}_2}.A_{{\h}_1}\nabla{v} -\nabla{u}.\AA \nabla{v} \right].
\end{equation}
where $\AA$ is defined in \eqref{definition:AA}.

%%%%%%%%%%%%%%%%%%%%%%%%%%%%%%
\subsection{Derivation of \eqref{derivee:seconde:etatD} from the weak formulation.}
\label{section:preuve:differentibilite:etape2}
%%%%%%%%%%%%%%%%%%%%%%%%%%%%%%

We want to make explicit the problem solved by  $(u')'$. To
achieve this, we should write the right hand side
\begin{equation*}
F = \int_\Omega \sigma \left[ \nabla{\dot{u}_1}.A_{{\h}_2}\nabla{v} + \nabla{\dot{u}_2}.A_{{\h}_1}\nabla{v} -\nabla{u}.\AA \nabla{v} \right],
\end{equation*}
as the sum of an integral with $\nabla v$ in factor and an integral of a divergence to identify the jump conditions on $\partial\omega$. To that end, we will use algebraic identities that involve second order derivatives of $u, \dot{u}_i$ and of the test function $v \in \mathcal{D}(\Omega)$. Using Lemma \ref{fait:technique0}, we obtain:
\begin{equation*}
\int_\Omega \sigma  \nabla{\dot{u}_1}.A_{{\h}_2}\nabla v   =  \int_\Omega \sigma \Big[ \nabla ({\h}_2.\nabla \dot{u}_1 ). \nabla v + \nabla ({\h}_2.\nabla v )\nabla \dot{u}_1 - \div{(\nabla \dot{u}_1 . \nabla v ) {\h}_2} \Big],
\end{equation*}
\begin{equation*}
 \int_\Omega \sigma  \nabla{\dot{u}_2}.A_{{\h}_1}\nabla{v}  =  \int_\Omega \sigma \Big[ \nabla ({\h}_1.\nabla \dot{u}_2 ). \nabla v + \nabla ({\h}_1.\nabla v )\nabla \dot{u}_2 - \div{(\nabla \dot{u}_2 . \nabla v ) {\h}_1} \Big].
\end{equation*}
Concerning the remaining terms, we use Lemma \ref{fait:technique1} to get
\begin{equation*}
\begin{split}
\int_{\Omega}\sigma \nabla{u}.\AA   \nabla v
=\int_{\Omega} & \sigma  ~\div{({\h}_2.\nabla u) A_{{\h}_1} \nabla v +({\h}_2.\nabla v) A_{{\h}_1} \nabla u -(A_{{\h}_1} \nabla u.\nabla v){\h}_2} \\
&~~- \sigma \Big[ ({\h}_2.\nabla u) \div{A_{{\h}_1} \nabla v} +  ({\h}_2.\nabla v) \div{A_{{\h}_1} \nabla u} \Big].
\end{split}
\end{equation*}
We apply Lemma \ref{fait:technique2} and gather the expressions
obtained for $F$. 
\begin{equation}
\begin{split}
F=\int_\Omega \sigma &\Big[ \nabla \left({\h}_1.\nabla \dot{u}_2 + {\h}_2.\nabla\dot{u}_1 \right) .\nabla v
+ \nabla ({\h}_2.\nabla v).\nabla \dot{u}_1 +\nabla ({\h}_1.\nabla v).\nabla \dot{u}_2 \Big] \\
& + \int_\Omega \sigma ~\div{ (A_{{\h}_1}\nabla u.\nabla v - \nabla \dot{u}_1. \nabla v) {\h}_2  - (\nabla \dot{u}_2. \nabla v) {\h}_1}\\
& + \int_\Omega \sigma \Big[ ({\h}_2.\nabla v) \Delta ({\h}_1. \nabla u) - \div{({\h}_2.\nabla v) A_{{\h}_1}\nabla u} - \nabla ({\h}_2.\nabla u). A_{{\h}_1}\nabla v  \Big].
\end{split}
\end{equation}
Using \eqref{formule:magique1}, we remove the dependency on
$A_{{\h}_1}\nabla v$:
\begin{equation*}
\nabla ({\h}_2.\nabla u). A_{{\h}_1}\nabla v  = \nabla ({\h}_1.\nabla ({\h}_2.\nabla u) ). \nabla v + \nabla ({\h}_1.\nabla v )\nabla ({\h}_2.\nabla u) - \div{(\nabla ({\h}_2.\nabla u) . \nabla v ) {\h}_1}.
\end{equation*}
Therefore, we write $F=F_1+F_2$ where
\begin{eqnarray}
F_1 & = &  \int_\Omega \sigma \Big[ \nabla \left({\h}_1.\nabla \dot{u}_2 + {\h}_2.\nabla\dot{u}_1\right) -\nabla ({\h}_1.\nabla ({\h}_2.\nabla u) )\Big] .\nabla v , \label{expression:F1} \\
F_2 & = & \int_\Omega \sigma \Big[ \nabla ({\h}_1.\nabla v).\nabla (\dot{u}_2 -{\h}_2.\nabla u) + \nabla ({\h}_2.\nabla v).\nabla \dot{u}_1 +  ({\h}_2.\nabla v) \Delta ({\h}_1. \nabla u)  \Big] \notag \\
&& + \int_\Omega \sigma ~\div{ (A_{{\h}_1}\nabla u.\nabla v-\nabla \dot{u}_1. \nabla v) {\h}_2    + \left(\nabla ({\h}_2.\nabla u) . \nabla v - \nabla \dot{u}_2. \nabla v\right) {\h}_1-({\h}_2.\nabla v) A_{{\h}_1}\nabla u} \notag.
\end{eqnarray}
The connection between second order material and shape derivatives is given by:
\begin{equation*}
\ddot{u}_{1} = (u'_1)'_2 + \h_1.\nabla \dot{u}_2 + \h_2.\nabla \dot{u}_1-\h_1.\nabla(\h_2.\nabla u),
%= (u'_1)'_2+ \h_1.\nabla u'_2 + \h_2.\nabla \dot{u}_1.
\end{equation*}
incorporating this expression in \eqref{expression:F1}, we rewrite \eqref{derivee:seconde:materielle:faible} as:
\begin{equation}
\label{reecriture:formualation:faible:derivee:forme}
\forall v \in H^1_0(\Omega),~~\int_{\Omega} \sigma \nabla (u'_1)'_2.\nabla v = F_2.
\end{equation}
Testing it against $v \in \mathcal{D}(\Omega\setminus\partial\omega)$, we get $\Delta (u'_1)'_2 =0$ in $\Omega\setminus\overline{\omega}$ and in $\omega$. We now deduce the jump conditions for $(u_1')'_2$. To obtain the jump of the potential, we simply write that $\ddot{u}_1 \in \sH^1_0(\Omega)$, hence $[\ddot{u}_1] =0$ on $\partial\omega$ and then
\begin{equation*}
 [(u'_1)'_2] = -\h_1.\nabla u_2' -\h_2.\nabla \dot{u}_1.
 \end{equation*}
To express the jump of the flux, we then apply the Gauss formula in \eqref{reecriture:formualation:faible:derivee:forme} to get
\begin{equation}\label{F2}
- \int_{\partial\omega} [\sigma \partial_{\n} (u'_1)'_2] v = F_2.
\end{equation}
The second term $F_2$ contains all the jumps of the flux  on the interface
$\partial\omega$.

 \textbf{A simplified  expression of $F_2$.}  To get a simplified formula for $F_2$ under a boundary integral, some lengthy but straightforward calculations are needed.
 We summarize the result by means of the following lemma
\begin{lemma}
One has:
\begin{equation}
\begin{split}
F_2 = \int_{\partial \omega} &\divt{ 2 h_{2,n} h_{1,n}  D\n \left[\sigma \nablat  u \right]
      - h_{2,n} \n. \nabla h_{1,n} \left[\sigma \nablat u \right]+  h_{2,n}\ct{{\h}_1} . D\n \,\n ~ \left[\sigma \nablat u \right]  }  v   \\
&+\int_{\partial \omega} \divt{\ct{{\h}_1}.\nablat(h_{2,n}) \left[\sigma \nablat u \right]-h_{1,n}h_{2,n} H\left[\sigma \nablat u \right]}v\\
&-\int_{\partial \omega} \left(\divt{h_{2,n}\left[\sigma \nablat u'_1 \right] }+\divt{h_{1,n}\left[\sigma \nablat u'_2 \right] }\right)v.
\end{split}
\end{equation}
\end{lemma}
\begin{proofof} {lemma} First, write :
\begin{equation*}
\begin{split}
\int_\Omega \sigma \nabla ({\h}_1.\nabla v).\nabla (\dot{u}_2 -{\h}_2.\nabla u)
&= \sigma_1 \int_{\Omega \setminus \overline{\omega}} \nabla ({\h}_1.\nabla v).\nabla u'_2
+ \sigma_2 \int_{\omega} \nabla ({\h}_1.\nabla v).\nabla u'_2\\
& =   - \int_{\partial \omega} [\sigma  \partial_{\n} u' _2] ({\h}_1.\nabla v)
\end{split}
\end{equation*}
Note that the normal vector is oriented from $\omega$ to $\Omega \setminus \overline{\omega}$.
In the same spirit, we write
\begin{equation*}
\nabla ({\h}_2.\nabla v).\nabla \dot{u}_1 +  ({\h}_2.\nabla v) \Delta ({\h}_1. \nabla u) = \nabla ({\h}_2.\nabla v).\nabla (\dot{u}_1- {\h}_1.\nabla u) + \div{ ({\h}_2.\nabla v). \nabla ({\h}_1.\nabla u)}.
\end{equation*}
By a argument of symmetry, we then can write:
\begin{equation*}
\int_\Omega \sigma \nabla ({\h}_2.\nabla v).\nabla (\dot{u}_1 -{\h}_1.\nabla u)=   - \int_{\partial \omega} [\sigma  \partial_{\n} u' _1] ({\h}_2.\nabla v).
\end{equation*}
To drop the dependency in $A_{{\h}_1}$, we use \eqref{formule:magique1} and get after expansion:
\begin{equation*}
\div{ \left(A_{{\h}_1} \nabla u . \nabla v\right) {\h}_2 }=\div{ \left( \nabla ({\h}_1.\nabla v ). \nabla u + \nabla ({\h}_1.\nabla u )\nabla v \right) {\h}_2  }
-\div{  \div{(\nabla u . \nabla v ) {\h}_1}  {\h}_2 };
\end{equation*}
\begin{eqnarray*}
\div{ ({\h}_2.\nabla v) A_{{\h}_1}\nabla u} &  = &\nabla  ({\h}_2.\nabla v). A_{{\h}_1}\nabla u +  ({\h}_2.\nabla v)\div{ A_{{\h}_1}\nabla u} \\
& =&  \nabla ({\h}_1.\nabla ({\h}_2.\nabla v)). \nabla u + \nabla ({\h}_1.\nabla u )\nabla ({\h}_2.\nabla v) +  ({\h}_2.\nabla v)\Delta({\h}_1.\nabla u)\\
& & ~~- \div{\left(\nabla ({\h}_2.\nabla v).  \nabla u \right) {\h}_1}\\
& = & \nabla \left({\h}_1.\nabla ({\h}_2.\nabla v)\right). \nabla u +  \div{ ({\h}_2.\nabla v)\nabla({\h}_1.\nabla u)-\left(\nabla ({\h}_2.\nabla v).  \nabla u \right) {\h}_1 }  .
\end{eqnarray*}
After integrating by parts, we  conclude thanks to  the state
equation  and obtain
\begin{equation*}
\int_\Omega \sigma \nabla \left({\h}_1.\nabla ({\h}_2.\nabla v)\right). \nabla u = -\int_\Omega \left({\h}_1.\nabla ({\h}_2.\nabla v)\right) \div{\sigma \nabla u} =0
\end{equation*}
We substitute the shape derivative $u'$ to the material one $\dot{u}$:
\begin{eqnarray*}
F_2 & = & -\int_{\partial \omega} [\sigma \partial_{\n} u'_1] ({\h}_2.\nabla v) + [\sigma  \partial_{\n} u'_2] ({\h}_1.\nabla v) - \int_{ \omega} \sigma ~\div{ \div{(\nabla u .\nabla v){\h}_1} {\h}_2}\\
&&  + \int_{ \omega} \sigma ~\div{ \left(\left(\nabla({\h}_1.\nabla v).\nabla u) {\h}_2 +(\nabla ({\h}_2.\nabla v).\nabla u \right) {\h}_1\right)-\left( (\nabla u'_2.\nabla v) {\h}_1 + ( \nabla u'_1 .\nabla v) {\h}_2\right)}.
 \end{eqnarray*}
First, we  use the continuity of the flux on $\partial\omega$,
then we integrate by parts on $\partial\omega$ and finally we
incorporate the expressions of the jumps of the shape derivatives
$u'$ to obtain
 \begin{eqnarray*}
\int_{\Omega}\sigma  \div{{\h}_1.\left(\nabla({\h}_2.\nabla v).\nabla u\right)} &=&-\int_{\partial \omega} \left[\sigma \nabla u.\nabla({\h}_2.\nabla v)\right] h_{1,n} =-\int_{\partial \omega} \left[\sigma \nablat u \right]h_{1,n} \nablat ({\h}_2.\nabla v) \\
&= &\int_{\partial \omega} \divt{\left[\sigma \nablat u\right]h_{1,n}} {\h}_2.\nabla v
= \int_{\partial \omega} \left[\sigma \partial_{\n} u'_1\right] {\h}_2.\nabla v.
\end{eqnarray*}
This leads to a simplified expression for $F_2$:
\begin{equation*}
F_2  =  - \int_{ \omega} \sigma ~\div{ \div{(\nabla u .\nabla v){\h}_1}\h_2 + \left( (\nabla u'_1 .\nabla v) {\h}_2 +  (\nabla u'_2.\nabla v) {\h}_1\right)  }  .
\end{equation*}
Let us study each term of this sum. Using  Gauss formula and
integrating by parts on the manifold $\partial\omega$, we obtain
\begin{equation*}
\begin{split}
\int_{ \omega} \sigma ~\div{ \nabla u'_1 .\nabla v) {\h}_2 }
   &=  - \int_{\partial \omega} h_{2,n} \left[\sigma \nabla u'_1.\nabla v\right]=  - \int_{\partial \omega} h_{2,n} \left[\sigma \partial_{\n} u'_1\right] \partial_{\n} v - \int_{\partial \omega} h_{2,n} \left[\sigma \nablat u'_1\right] \nablat v \\
&=  -  \int_{\partial \omega} h_{2,n} \left[\sigma \partial_{\n} u'_1\right] \partial_{\n} v+\int_{\partial \omega} \divt{ h_{2,n} \left[\sigma \nablat u'_1\right]}  v.
\end{split}
\end{equation*}
By symmetry, we also get:
\begin{equation*}
\int_{ \omega} \sigma ~\div{ \nabla u'_2 .\nabla v) {\h}_1 } =
-  \int_{\partial \omega} h_{1,n} \left[\sigma \partial_{\n} u'_2 \right] \partial_{\n} v+\int_{\partial \omega} \divt{ h_{1,n} \left[\sigma \nablat u'_2 \right]}  v
\end{equation*}
We now turn to the term with a double divergence. We first write
it as a boundary integral thanks to Gauss formula as
\begin{equation*}
\int_{ \omega} \sigma ~\div{ \div{(\nabla u .\nabla v){\h}_1} {\h}_2} = \int_{\partial \omega} h_{2,n} \div{{\h}_1\left[\sigma (\nabla u .\nabla v)\right]},
\end{equation*}
then, we use \eqref{definition:divergence:tangentielle} to introduce the tangential operators
\begin{equation*}
  \int_{ \omega} \sigma ~\div{ \div{(\nabla u .\nabla v){\h}_1} {\h}_2}=  \int_{\partial \omega} h_{2,n} \divt{{\h}_1\left[\sigma (\nabla u .\nabla v)\right]}+\int_{\partial \omega} h_{2,n} D({\h}_1\left[\sigma (\nabla u .\nabla v)\right])\n.\n.
  \end{equation*}
We study each of these terms. We start with the one involving
tangential derivatives: we expand the tangential divergence to
incorporate the jump relation for the state $u$.
\begin{eqnarray*}
 \divt{{\h}_1\left[\sigma (\nabla u .\nabla v)\right]}
 & = & \divt{{\h}_1} \left[\sigma (\nabla u .\nabla v)\right] + {\h}_1. \nablat \left[\sigma \nabla u. \nabla v\right]\\
 & =  & \divt{{\h}_1} \left[\sigma \nablat u\right]. \nablat v + {\h}_1. \nablat \left[\sigma \nablat u. \nablat v\right].
\end{eqnarray*}
Then, the first term  becomes:
\begin{equation*}
 \int_{\partial \omega} h_{2,n} \divt{{\h}_1\left[\sigma (\nabla u .\nabla v)\right]} = \int_{\partial \omega}h_{2,n} \divt{{\h}_1} \left[ \sigma \nablat u \right] . \nablat v + \int_{\partial \omega} h_{2,n}  {\h}_1. \nablat \left[\sigma \nablat u \nablat{v}\right].
\end{equation*}
We use the integration by parts formula \eqref{formule:integration:par:parties} to get:
\begin{equation*}
\begin{split}
 \int_{\partial \omega} &h_{2,n} \divt{{\h}_1\left[\sigma (\nabla u .\nabla v)\right]} \\
& = \int_{\partial \omega}  h_{1,n} h_{2,n}H \left[ \sigma \nablat u \right]. \nablat v   - \divt{ \divt{{\h}_1} h_{2,n} \left[ \sigma \nablat u \right]} v - \divt{{\h}_1 h_{2,n}}  \left[ \sigma \nablat u \right]. \nablat{v}  \\
& =\int_{\partial \omega}  \divt{ \left( \divt{ {\h}_1 h_{2,n} } -\divt{{\h}_1} h_{2,n} - h_{1,n} h_{2,n} H\right) \left[\sigma \nablat u \right]}  v.
\end{split}
\end{equation*}
Expanding
\begin{eqnarray*}
  \divt{\divt{{\h}_1 h_{2,n}} \left[\sigma \nablat u \right]}v &=&  \divt{\divt{{\h}_1}h_{2,n} \left[\sigma \nablat u \right]+ {\h}_1.\nablat(h_{2,n}) \left[\sigma \nablat u \right]}v\\
&=&  \divt{\divt{{\h}_1}h_{2,n} \left[\sigma \nablat u \right]}v+ \divt{\ct{{\h}_1}\nablat h_{2,n} \left[\sigma \nablat u \right]}v,
\end{eqnarray*}
we obtain the new expression:
\begin{equation}
 \int_{\partial \omega} h_{2,n} \divt{{\h}_1\left[\sigma (\nabla u .\nabla v)\right]} =
 \int_{\partial \omega}  \divt{ \left( \ct{{\h}_1}\nablat h_{2,n} - h_{1,n} h_{2,n} H\right) \left[\sigma \nablat u \right]} v.
\end{equation}
Now, we consider the term involving normal components.  We
have
\begin{equation}
\label{formule:magique3}
\begin{split}
\n.D({\h}_1\left[\sigma \nabla(\nabla u .\nabla v)\right])\n
&= \n.\nabla (h_{1,n}\left[\sigma \nabla u .\nabla v\right]) -  \left[\sigma \nabla u .\nabla v\right] \ct{{\h}_1}.D\n \,\n \\
&= \n. \nabla(h_{1,n}) \left[\sigma \nablat u \right] \nablat v+ h_{1,n}\n. \nabla( \left[\sigma \nabla u .\nabla v \right]).
\end{split}
\end{equation}
Then, we get
\begin{equation*}
\begin{split}
\int_{\partial \omega}& h_{2,n} D({\h}_1\left[\sigma (\nabla u .\nabla v)\right])\n.\n
 =  \int_{\partial \omega} h_{2,n}  \n. \nabla(h_{1,n})  \left[\sigma \nablat u \right].\nablat v   + h_{2,n} h_{1,n}\n. \nabla( \left[\sigma \nabla u .\nabla v \right])\\
&=   \int_{\partial \omega} - \divt{h_{2,n} \n. \nabla(h_{1,n}) \left[\sigma \nablat u \right]} v+h_{2,n} h_{1,n}\n. \nabla( \left[\sigma \nabla u .\nabla v \right]).
\end{split}
\end{equation*}
A straightforward  calculus leads to
\begin{equation*}
\begin{split}
\n. \nabla( \left[\sigma \nabla u .\nabla v \right])&= \n.\left(\left[\sigma D^2u \nabla v\right]+ D^2v \left[\sigma \nabla u \right]  \right)\\
&=\n.\left( \partial_{\n} v \left[\sigma D^2u \right] \n + \left[\sigma D^2u\right] \nablat v  + D^2 v   \left[\sigma \nablat u \right] \right)\\
&= \partial_{\n} v  \left[\sigma \cfrac{\partial^2 u}{\partial n^2}\right]+ \n .  \left[\sigma D^2 u \right]\nablat v + \n. D^2 v \left[\sigma \nablat u \right].
\end{split}
\end{equation*}
where $D^2 u$ is the Hessian matrix of $u$. From \eqref{lien:laplacien:laplace-beltrami} and from
the jump conditions for the state $u$, we deduce  that
\begin{equation*}
 \left[\sigma \cfrac{\partial^2 u}{\partial n^2}\right]=- \left[\sigma \Deltat u \right].
\end{equation*}
When one differentiates the relation expressing  the continuity of the flux for the state along the tangential direction
$\nablat v$, one gets (\cite{HenrotPierre}, p 235):
\begin{equation*}
0 =  \nabla [\sigma \partial_{\n} u].\nablat v = [\sigma D^2 u ]\nablat v.\n + [\sigma \nabla u]. (D\n \,\nablat v).
\end{equation*}
In the same spirit, it comes that
\begin{equation}
\label{derivation:suivant:tangente}
\nabla \partial_{\n} v . [\sigma \nablat u] = D^2 v [\sigma \nablat u].\n + \nabla v. (D\n [\sigma \nablat u]).
\end{equation}
Since $D\n$ is a symmetric matrix and $D\n\,\n=0$, one checks $\nabla v. (D\n [\sigma \nablat u])= [\sigma \nabla u]. (D\n \,\nablat v) $. Then
\begin{equation*}
\n. \nabla( \left[\sigma \nabla u .\nabla v \right]))
= - \left[\sigma \Deltat u \right]\partial_{\n} v
- 2 D\n \left[\sigma \nablat  u \right] . \nablat  v
+\left[\sigma \nablat  u \right] \nablat \partial_{\n} v
\end{equation*}
We integrate this expression on $\partial \omega$ and obtain after some
integration by parts:
\begin{equation*}
\begin{split}
\int_{\partial \omega}& h_{2,n} h_{1,n}\n. \nabla( \left[\sigma \nabla u .\nabla v \right]) \\
& = -\int_{\partial \omega}h_{2,n} h_{1,n} \left[\sigma \Deltat u \right]\partial_{\n} v
  +\int_{\partial \omega}h_{2,n} h_{1,n} \left[\sigma \nablat  u \right] \nablat \partial_{\n} v  -2 \int_{\partial \omega}h_{2,n} h_{1,n} D\n \left[\sigma \nablat  u \right] . \nablat  v ,\\
& = - \int_{\partial \omega} \left[ h_{2,n} h_{1,n} \left[\sigma \Deltat u \right]
  + \divt{h_{2,n} h_{1,n} \left[\sigma \nablat  u \right] } \right] \partial_{\n} v  +2  \int_{\partial \omega} \divt  {h_{2,n} h_{1,n}  D\n \left[\sigma \nablat  u \right] }  v  .
\end{split}
\end{equation*}
Hence
\begin{equation*}
\begin{split}
\int_{\partial \omega} h_{2,n}& D({\h}_1\left[\sigma (\nabla u .\nabla v)\right])\n.\n
 = - \int_{\partial \omega} \left[ h_{2,n} h_{1,n} \left[\sigma \Deltat u \right]
  + \divt{h_{2,n} h_{1,n} \left[\sigma \nablat  u \right] } \right] \partial_{\n} v \\
& +
   \int_{\partial \omega} \divt{ 2 h_{2,n} h_{1,n}  D\n \left[\sigma \nablat  u \right]
      - h_{2,n} \n. \nablat h_{1,n} \left[\sigma \nablat u \right]   }  v   .
\end{split}
\end{equation*}
\begin{equation*}
\begin{split}
 \int_{ \omega} \sigma& \div{ \div{(\nabla u .\nabla v){\h}_1} {\h}_2}=  \int_{\partial \omega} \left[ h_{2,n} h_{1,n} \left[\sigma \Deltat u \right]
  + \divt{h_{2,n} h_{1,n} \left[\sigma \nablat  u \right] } \right] \partial_{\n} v \\
& -
   \int_{\partial \omega} \divt{ 2 h_{2,n} h_{1,n}  D\n \left[\sigma \nablat  u \right]
      - h_{2,n} \n. \nabla h_{1,n} \left[\sigma \nablat u \right]  }  v   \\
&-\int_{\partial \omega} \divt{\ct{{\h}_1}.\nablat(h_{2,n}) \left[\sigma \nablat u \right]-h_{1,n}h_{2,n} H\left[\sigma \nablat u \right]}v.
\end{split}
\end{equation*}
Gathering all the terms, we write $F_2$ as:
\begin{equation*}
\begin{split}
F_2= \int_{\partial \omega} &\divt{ 2 h_{2,n} h_{1,n}  D\n \left[\sigma \nablat  u \right]
      +\left( \ct{{\h}_1}.\nablat(h_{2,n})  - h_{2,n} \n. \nabla h_{1,n}  -h_{1,n}h_{2,n} H \right) \left[\sigma \nablat u \right] }  v   \\
&-\int_{\partial \omega} \left(\divt{h_{2,n}\left[\sigma \nablat u'_1 \right] }+\divt{h_{1,n}\left[\sigma \nablat u'_2 \right] }\right)\\
&-\int_{\partial \omega} \left( h_{2,n} h_{1,n} \left[\sigma \Deltat u \right]
  + \divt{h_{2,n} h_{1,n} \left[\sigma \nablat  u \right] } \right) \partial_{\n} v \\
 & -\int_{\partial \omega} \left( h_{1,n} \divt{h_{2,n}\left[\sigma \nablat u \right]} + h_{2,n}\divt{h_{1,n}\left[\sigma \nablat u \right]}\right)\partial_{\n} v.
\end{split}
\end{equation*}
We end the proof after expanding  the tangential divergence of the last term of $F_2$.
\end{proofof}

Let us return to the weak formulation \eqref{F2} of the derivative.  By identification, we get
\begin{equation*}
\begin{split}
 [\sigma \partial_{\n} (u'_1)'_2]&=  \divt{h_{2,n}\left[\sigma \nablat u'_1 \right] }+\divt{h_{1,n}\left[\sigma \nablat u'_2 \right] } - \divt{ h_{2,n} h_{1,n}  (2 D\n -H I) \left[\sigma \nablat  u \right] } \\
  &-  \divt{\ct{{\h}_1}.\nablat(h_{2,n}) \left[\sigma \nablat u \right]      - h_{2,n} \n. \nabla h_{1,n} \left[\sigma \nablat u \right]+  h_{2,n}\ct{{\h}_1} . D\n \,\n ~ \left[\sigma \nablat u \right]  } .
\end{split}
\end{equation*}
It remains to  compute  the jump of the flux for the second order derivative. Since
\begin{equation}\label{relation:potentiel:derivee:seconde}
 u''_{1,2} = (u'_1)'_2 - u'_{D\h_1\,\h_2}
\end{equation}
where $u'_{D\h_1\,\h_2}$ is the first shape derivative of $u$ in the direction of the vector field $ D\h_1\,\h_2$. Thanks  to \eqref{derivee:etatD}, we can write the jump  under the form
\begin{equation}\label{relation:flux:derivee:seconde}
 [\sigma \partial_{\n} u''_{1,2}] =  [\sigma \partial_{\n}(u'_1)'_2 ]-
 [\sigma \partial_{\n}u'_{D\h_1\,\h_2}] =[\sigma \partial_{\n}(u'_1)'_2 ] -
 \divt{ D\h_1 \h_2 .\n [\sigma \nablat u]}.
\end{equation}
Let us  split the field $\h_2$ in two parts: $D\h_1\,\h_2 .\n =
h_{2,n} \n.D\h_1\,\n +  D\h_1 \, \ct{h_2}.\n $. In the spirit of
\eqref{derivation:suivant:tangente}, we obtain
\begin{equation}\label{Jacobien:suivant:tangente}
 D\h_1 \ct{h_2}.\n = \nablat h_{1,n}.\ct{h_2} - \ct{h_1}. D\n \, \ct{h_2}.
\end{equation}
Thanks to \eqref{formule:magique3}, the jump $ [\sigma \partial_{\n}u'_{D\h_1\,\h_2}]$ then can be written under the form
\begin{equation*}
[\sigma \partial_{\n}u'_{D\h_1\,\h_2}] = \divt{ ( h_{2,n} \n.\nabla h_{1,n}  +  \nablat h_{1,n}.\ct{h_2} -
\ct{h_1}. D\n \, \ct{h_2})[\sigma \nablat u]} .
\end{equation*}
Gathering all the terms, simplifications occur and we get:
\begin{equation*}
\begin{split}
 [\sigma \partial_{\n} u''_{1,2}]=&  \divt{h_{2,n}\left[\sigma \nablat u'_1 \right] + h_{1,n}\left[\sigma \nablat u'_2 \right] }-  \divt{(\ct{{\h}_1}.\nablat h_{2,n} +  \nablat h_{1,n}.\ct{h_2})  \left[\sigma \nablat u \right]   }\\
 & - \divt{ h_{2,n} h_{1,n}  (2  D\n-H I) \left[\sigma \nablat  u \right] }  + \divt{\ct{h_1}. D\n \, \ct{h_2})[\sigma \nablat u]}.
\end{split}
\end{equation*}
To get the jumps of the potential, we  use  \eqref{relation:potentiel:derivee:seconde} and obtain
\begin{equation*}
\begin{split}
\left[u''_{1,2}\right] &= \left[(u'_1)'_2\right] -
\left[u'_{D\h_1\,\h_2}\right]=-\h_1.\left[\nabla
u'_{2}\right]-\h_2.\left[\nabla
\dot{u_1}\right]-\left[u'_{D\h_1\,\h_2}\right]\\
&=-h_{2,n}\left[\partial_\n u'_{1}\right] -h_{1,n}\left[\partial_\n
u'_2 \right] -\ct{h_2}.\left[\nablat
u'_{1}\right]-\ct{h_1}.\left[\nablat u'_2\right]\\
&~~~~~~-h_{2,n} \n.\left[ \nabla (\h_1.\nabla u)\right]-h_{1,n}
\n.\left[ \nabla (\h_2.\nabla u)\right]-
\left[u'_{D\h_1\,\h_2}\right]
\end{split}
\end{equation*}
Thanks to   the jump of the potential for  the first order shape
derivative given in \eqref{derivee:etatD}, it comes that
\begin{equation*}
\ct{h_2}.\left[\nablat u'_1\right]=-\ct{h_2}.\left[\nabla(\h_1.\nabla u)\right]~~~~\text{and}~~~~\ct{h_1}.\left[\nablat u'_2\right]=-\ct{h_1}.\left[\nabla(\h_2.\nabla u)\right]
\end{equation*}
and then:
\begin{equation}\label{saut:derivee:seconde:potentiel}
\left[u''_{1,2}\right] =-h_{2,n}\left[\partial_\n u'_{1}\right]
-h_{1,n}\left[\partial_\n u'_2\right] -h_{2,n} \n.\left[ \nabla
(\h_1.\nabla u)\right]+\ct{h_1}.\left[\nabla(\h_2.\nabla u)\right]
- \left[u'_{D\h_1\,\h_2}\right]
\end{equation}
Computing the other  jumps that  appeared in the former
expression, we get
\begin{eqnarray*}
\left[ \nabla(\h_2.\nabla u)\right] & =& (D\h_2)^T\left[
\nabla u\right]+\left[ D^2 u\right]\h_2.\\
\ct{\h_1}.\left[\nabla(\h_2.\nabla u) \right]& = & \n.D\h_2\ct{h_1}\left[ \partial_\n u\right] +h_{2,n} \ct{h_1}.\left[D^2 u\right]\n+\ct{h_1}.\left[D^2 u\right]\ct{h_2}.\\
h_{2,n} \n.\left[\nabla(\h_1.\nabla u) \right]& = & h_{2,n}\left[ \partial_\n u\right]  \n.D\h_1\n +h_{2,n} h_{1,n}\n.\left[D^2 u\right]\n+h_{2,n}\n.\left[D^2 u\right]\ct{h_1}.\\
\left[ u'_{D\h_1\, \h_2}\right]& = & -D\h_1\,\h_2 .\n \left[ \partial_\n u
 \right]  =-\left(h_{2,n} \n.D\h_1\n+\n.D\h_1\ct{h_2}\right) \left[ \partial_\n u\right]
 \end{eqnarray*}
 With the help  of formula  \eqref{Jacobien:suivant:tangente}, we obtain:
\begin{equation*}
\begin{split}
-h_{2,n} \n.\left[ \nabla (\h_1.\nabla u)\right]+\ct{h_1}.&\left[\nabla(\h_2.\nabla u)\right] - \left[u'_{D\h_1\,\h_2}\right]
= \left(\nablat h_{1,n}.\ct{h_2}+\nablat h_{2,n}.\ct{h_1}\right)\left[\partial_\n u\right]\\
& -2 \ct{h_1}.D\n\ct{h_2} \left[\partial_\n u \right] +\ct{h_1}.\left[D^2 u \right]\ct{h_2}-h_{2,n} h_{1,n} \n.\left[D^2 u\right]\n.
\end{split}
\end{equation*}
\begin{equation*}
\n.\left[D^2 u\right]\n=-\left[ \Deltat u\right]-H \left[\partial_\n u\right]=-H \left[\partial_\n u\right],
\end{equation*}
\begin{equation*}
\ct{h_1}.\left[D^2 u\right]\ct{h_2}=\ct{h_1}.D(\left[\nabla u \right])\ct{h_2}
 =\ct{h_1}.D(\left[\partial_\n u \right]\n)\ct{h_2}=\left[\partial_\n u \right]\ct{h_1}.D\n\ct{h_2}.
\end{equation*}
Finally, we gather the results of these computations to write
\begin{equation}
\begin{split}
\left[u''_{1,2}\right] = -\left(h_{2,n}\left[\partial_\n u'_{1}\right] +h_{1,n}\left[\partial_\n u'_2\right]\right)&
     + \left(\nablat h_{1,n}.\ct{h_2}+\nablat h_{2,n}.\ct{h_1}\right)\left[\partial_\n u \right]\\
&+\left(h_{2,n} h_{1,n} H -\ct{h_1}.D\n\, \ct{h_2}\right)\left[\partial_\n u\right])
\end{split}
\end{equation}

%%%%%%%%%%%%%%%%%%%%%%%%%%%%%
\subsection{How to recover \eqref{derivee:seconde:etatD} by formal differentiation of the boundary conditions.}
\label{section:preuve:differentibilite:etape3}
%%%%%%%%%%%%%%%%%%%%%%%%%%%%%

The aim of this section is to retrieve the expression of the flux jump  $[\sigma \partial_\n u'']$ by computing the normal  derivatives of each of the expressions
$\dot{\overline{[\sigma \nabla u'].\n}}$ and $\dot{\overline{\divt{h_{1,n}
[\sigma \nablat u]}}}$. Since
\begin{equation*}
[\sigma\nabla u'].\n = \divt{h_{1,n} [\sigma\nablat u]} = h_{1,n} [\sigma \Deltat u] + \nablat h_{1,n} . [\sigma \nablat u],
\end{equation*}
 then, we get
\begin{equation}
\label{derivee:meterielle:saut:de:flux}
\dot{\overline{[\sigma \nabla u'].\n}} = \dot{\overline{ h_{1,n}}} [\sigma \Deltat u] +  h_{1,n} \dot{\overline{[\sigma \Deltat u] }} +   \dot{\overline{\nablat h_{1,n}}} . [\sigma \nablat u]+ \nablat h_{1,n} .  \dot{\overline{[\sigma \nablat u]}}.
\end{equation}
In order  to avoid lengthy computations, we shall concentrate on each
normal  derivative appearing in the above formula. Some of the results are straightforward and
their proof will be left to the reader.
Combining  propositions  \eqref{derivative:normal} and  \eqref{derivee:materielle:gradient:tangentiel}, we conclude that
\begin{equation*}
 \dot{\overline{\nablat h_{1,n} }} = - \nablat ( \h_1.\nablat h_{2,n} ) + (D^2 h_{1,n} .\h_2)_\tau - \nabla h_{1,n} .\dot{\n}~ \n-\nabla h_{1,n}.\n~ \dot{\n}.
\end{equation*}
In the same manner, we also get
\begin{equation*}
 \dot{\overline{ [\sigma \nablat u]}} = [ \sigma \nablat u'_2]  +  ([\sigma D^2 u] .\h_2)_\tau - [\sigma \nablat u]  .\dot{\n}~\n-[\sigma \nablat u]\n~ \dot{\n}.
\end{equation*}
Hence, we can write
\begin{equation*}
 \dot{\overline{\h_1.\n}} = \h_2 .\nabla h_{,n} - \nablat h_{2,n} . \ct{\h_1}.
\end{equation*}
 It remains to simplify the terms $ A= (D^2 u .\h_2)_\tau  .\nablat h_{1,n} $ and
 $ B=[\sigma \nablat u]. (D^2 h_{1,n}.\h_2)_\tau$. We obtain:
 \begin{eqnarray*}
A&=& -[\sigma \nablat u]. (D\n\,  \nablat h_{1,n}) h_{2,n} + [\sigma \Deltat u] \nablat h_{1,n}.\ct{\h_2},\\
B&=& (D^2 h_{2,n}.\ct{\h_2} ).[\sigma \nablat u] + \nablat(\partial_\n h_{1,n} ). [\sigma\nablat u] h_{2,n}
 - [\sigma \nablat u].(D\n\, \nablat h_{1,n}) h_{2,n}.
\end{eqnarray*}
We tackle the computation of $(\partial_\n u')'$. For the sake of
clearness, we subdivide the work in several steps.

\emph{First step.} We compute $\dot{\overline{\divt{h_{1,n} [\sigma \nablat u]} }}$. We expand:
\begin{eqnarray*}
\dot{\overline{\divt{h_{1,n} [\sigma \nablat u]} }} & =& \dot{\overline{h_{1,n} [\sigma \Deltat u] }} + \dot{\overline{\nablat h_{1,n} . [\sigma \nablat u] }},\\
&=& \dot{\overline{h_{1,n}}} [\sigma \Deltat u] + h_{1,n} \dot{\overline{[\sigma \Deltat u] }} + \dot{\overline{\nablat h_{1,n}}}.[\sigma \nablat u] + \nablat h_{1,n}. \dot{\overline{[\sigma \nablat u]}}.
\end{eqnarray*}

Hence, after substitution, one gets
\begin{eqnarray}
\dot{\overline{ \divt{h_{1,n}  [ \sigma \nablat u] } }} & = & \divt{ h_{1,n} [ \sigma \nablat u'_2 ] + ( h_{2,n} \partial_\n h_{1,n} - \nablat h_{2,n}. \ct{\h_1}) [\sigma \nablat u] }  \nonumber\\
&&~~ + 2 [\sigma \Deltat u] \nablat h_{1,n}.\ct{\h_2}- \partial_\n h_{1,} [\sigma\nablat u].(D\n\, \ct{\h_2}) + [\sigma \nablat u] .(D^2 h_{1,n} .\ct{\h_2}) \nonumber \\
&& ~~ - 2 h_{2,n} [\sigma \nablat u] .(D\n\, \nablat h_{1,n}) + h_{1,n} \left( \dot{\overline{\sigma \Deltat u}} - [\sigma \Deltat u'_2] \right).
\label{etoile1}
\end{eqnarray}

\emph{Second step.} We compute $\dot{\overline{[\sigma \partial_\n
u'_1]}} $.  From the expression of $\dot{\n}$, we get after some
straightforward computations:
\begin{equation}
\label{etoile2}
\dot{\overline{[\sigma \partial_\n u'_1]}} = [\sigma \partial_\n (u'_1)'_2] + ([\sigma D^2 u'_1]\,\h_2).\n + [\sigma \nablat u'_1]. (D\n\, \ct{\h_2}-\nablat h_{2,n} ).
\end{equation}

\emph{Third step.} We compute $\sigma \partial_\n (u'_1)'_2$. From the jump condition on the flux of the
derivative \eqref{derivee:etatD} and \eqref{etoile1} and \eqref{etoile2}, we obtain:
\begin{eqnarray*}
[\sigma \partial (u'_1)'_2]&=& \divt{h_{1,n} [\sigma\nablat u'_2] + \left(h_{2,n} \partial_\n h_{1,n}- \nablat h_{2,n}.\ct{\h_1} \right) [\sigma \nablat u]   }+ 2\nablat h_{1,n}.\ct{\h_2} [\sigma \Deltat u] \\
&&-([\sigma D^2 u'_1]\,\h_2).\n + [\sigma \nablat u'_1]. \left( \nablat h_{2,n} - D\n\, \ct{\h_2} \right) - \partial_\n h_{1,n} [\sigma\nablat u].(D_n\,\ct{\h_2}) \\
&&+ (D^2 h_{1,n}\,\ct{\h_2}).[\sigma\nablat u]- h_{2,n} (D\n\, [\sigma\nablat u])£.\nablat h_{1,n}.
\end{eqnarray*}
Taking account of the following calculation,
\begin{equation*}
\begin{split}
-([\sigma D^2 u'_1]\,\h_2).\n + [\sigma \nablat u'_1] .& \nablat {h_2,n}  = - \left( h_{2,n} [\sigma D^2u'_1]\,\n+ [\sigma D^2u'_1]\,\ct{\h_2} \right).\n + [\sigma \nablat u'_1] . \nablat {h_2,n},\\
&= h_{2,n} \left([\sigma \Deltat u'_1] + H [\sigma \partial_\n u'_1] \right) +  [\sigma u'_1].\nablat h_{2,n} - ([\sigma D^2 u'1]\,\ct{\h_2}).\n,\\
&= \divt{h_{2,n} [\sigma \nablat u'_1]} + H h_{2,n} [\sigma \partial_\n u'_1] - ([\sigma D^2 u'_1]\,\ct{\h_2}).\n;
\end{split}
\end{equation*}
it comes
\begin{eqnarray}
[\sigma \partial_\n(u'_1)'_2] &=& \divt{h_{1,n} [\sigma\nablat u'_2] + h_{2,n} [\sigma\nablat u'_1]+  \left(h_{2,n} \partial_\n h_{1,n}- \nablat h_{2,n}.\ct{\h_1} \right) [\sigma \nablat u] } \nonumber \\
&& + 2 [\sigma \Deltat u] \nablat h_{1,n}.\ct{\h_2}+ H h_{2,n} [\sigma \partial_\n u'_1] -
([\sigma D^2 u'_1] \, \ct{\h_2}).\n \nonumber\\
&&-\left( [\sigma \nablat u'_1]+ \partial_\n h_{1,n} [\sigma\nablat u ] \right).(D\n\, \ct{\h_2}) + (D^2h_{1,n} \, \ct{\h_2} ) .[\sigma \nablat u] \nonumber\\
&& - 2 h_{2,n} \nablat h_{1,n}.(D\n\, [\sigma \nablat u ]) + h_{1,n}\left( \dot{\overline{[\sigma \Deltat u]}} -[\sigma \Deltat u'_2] \right) \label{etoile3}.
\end{eqnarray}
This formula remains  hard to handle. To get a more
convenient  one, we decide to derive tangentially to the direction
$\h_2$ the boundary identity
\begin{equation*}
[\sigma \partial_\n u'_1] = h_{1,n} [\sigma \Deltat u] + \nablat h_{1,n}.[\sigma \nablat u].
\end{equation*}
This leads to:
 \begin{equation}
 \label{etoile4}
 \begin{split}
([\sigma D^2 u'_1]\, \ct{\h_2}) .\n + &(D\n\, \ct{\h_2}).[\sigma \nablat u'_1] = \nablat h_{1,n}.\ct{\h_2} [\sigma \Deltat u] + h_{1,n} \nablat [\sigma\Deltat u].\ct{\h_2}\\
& + (D^2h_{1,n}\,\ct{\h_2}).[\sigma \nablat u] -\partial_\n h_{1,n} [\sigma \nablat u]. (D\n\, \ct{\h_2})+[\sigma \Deltat u ] \ct{\h_2}.\nablat h_{1,n}.
\end{split}
\end{equation}
From \eqref{formule:derivee:materielle:laplace:beltrami} and subtracting \eqref{etoile4} from \eqref{etoile3}, we can write
\begin{eqnarray*}
[\sigma \partial_\n (u'_1)'_2] &=& \divt{ h_{1,n} [\sigma\nablat u'_2] + h_{2,n} [\sigma\nablat u'_1]+  \left(h_{2,n} \partial_\n h_{1,n}- \nablat h_{2,n}.\ct{\h_1} \right)[\sigma \nablat u] } \\
&&+ \divt{h_{1,n} h_{2,n}(H I- 2D \n). [\sigma \nablat u] } - h_{1,n}\left( \nablat [\sigma \Deltat u].\ct{\h_2} + \Deltat [\sigma \nablat u].\ct{\h_2} \right)\\
&& +h_{1,n}\left(  \nablat \divt{\ct{\h_2}}.[\sigma \nablat u] -\divt{\left(\left( D{\h_2}+( D{\h_2})^T \right)[\sigma\nablat u] \right)_{\tau}}\right).
\end{eqnarray*}
From \eqref{formule:derivee:materielle:laplace:beltrami}, we obtain
\begin{equation}
\label{formule:derivee:materielle:laplace:beltrami:2}
\begin{split}
\dot{\overline{[\sigma \Deltat u ]}} = [\sigma\Deltat \dot{u}] &+ \nablat \divt{\ct{\h_2}}.[\sigma \nablat u ] + \nablat(H \h_{2,n}).[\sigma \nablat u ]\\&  - \divt{\left(\left( D{\h_2}+( D{\h_2})^T \right)[\sigma\nablat u]\right)_{\tau}},
\end{split}
\end{equation}
and using  the relation  between the material  and shape derivative, we get
\begin{equation*}
\dot{\overline{[\sigma \Deltat u ]}}= [\sigma\Deltat u']+ \nabla \left([\sigma \Deltat u]\right).\h_2 ~~\text{and}~~ [\sigma\Deltat \dot{u}]= [\sigma\Deltat u']+\Deltat\left([\sigma \nabla u].\h_2\right).
\end{equation*}
Injecting  these relations in \eqref{formule:derivee:materielle:laplace:beltrami:2} and  applying them for $\ct{\h_2}$, we get
\begin{equation*}
%\begin{split}
 \Deltat ([\sigma \nablat u]. \ct{\h_2} ) + \nablat \divt{\ct{\h_2}}.[\sigma\nablat u] = \nablat [\sigma \Deltat u]. \ct{\h_2}+ \divt{\left(\left( D{\h_2}+( D{\h_2})^T \right)[\sigma\nablat u] \right)_{\tau}}.
%\end{split}
\end{equation*}
This last fact allows us to conclude.
%============================================
%
%===========================================
\subsection{Justification of the formal computations.}
\label{section:preuve:differentibilite:etape4}

%Before starting to prove Theorem \ref{lem:derivee:etatD} and
%Theorem \ref{theoreme:derivee:seconde:etatD},
We have to justify rigorously that the right-hand sides of
\eqref{derivee:etatD},\eqref{derivee:etatN},\eqref{derivee:seconde:etatD}
make sense. They involve tangential derivatives of $\un$
and $\ud$ along the interface $\partial\omega$ up to the order
three. The existence of these derivatives is not clear \emph{a
priori} since the gradient of the solution has a discontinuity along
this interface.  Our first aim is to  precise the tangential
regularity along the interface $\partial \omega$ of the solution
$u$ of \eqref{equation:etat} with either Dirichlet or Neumann boundary conditions.

We should access to the trace of  $u$ on the interface
$\partial\omega$. Any numerical discretization needs also to
compute the state, its derivatives with respect to the shape and
the normal derivatives along the interface $\partial\omega$.  To
that end,  we introduce for any $\alpha \in H^{1/2}(\partial
\omega)$ and $\beta \in H^{-1/2}(\partial\omega)$  the following
boundary value problems
\begin{equation}
\label{probleme:transmission}
(D)\left\{
\begin{array}{rcl}
\Delta v  & =  & 0 \text{ in }\Omega\setminus\overline{\omega} \text{ and  in } \omega,   \\
  \left[v\right]& =  & \alpha \text{ on }\partial \omega, \\
  \left[\sigma \partial_{n} v\right]&  = & \beta  \text{ on }\partial \omega,\\
  v &=& f_1  \text{ on } \partial\Omega.
\end{array}
\right.
\text{ and }(N)
%\label{probleme:transmission:modele:Neumann}
\left\{
\begin{array}{rcl}
\Delta v  & =  & 0 \text{ in }\Omega\setminus\overline{\omega} \text{ and  in } \omega,   \\
  \left[v\right]& =  & \alpha \text{ on }\partial \omega, \\
  \left[\sigma \partial_{n} v\right]&  = & \beta  \text{ on }\partial \omega,\\
  \partial_\n v &=& g_1  \text{ on } \partial\Omega,
\end{array}
\right.
\end{equation}
where $(f_1,g_1)\in H^{1/2}(\partial\Omega)\times
H^{-1/2}(\partial\Omega)$. Note that for $\alpha=0$, $\beta=0$ and
$(f_1,g_1)=(f,g)$ then $(\ud)$ and $\un$ solve respectively (D)
and (N); furthermore the choice of
\begin{equation}
\label{choix:sauts}
\alpha=\cfrac{\left[\sigma\right]}{\sigma_2} h_n \partial_\n u^+
\text{ and }\beta=\left[\sigma\right] \divt{h_n \nablat u} 
\end{equation}
leads
to \eqref{derivee:etatD} and \eqref{derivee:etatN} when we take
$(f_1,g_)=(0,0).$

\textbf{Existence of solutions to (D) and (N).} To study these problems, we use
the integral representation in terms of layer potentials. In a
first step, we recall some definitions. The Newtonian potential
$\Gamma$  is defined as:
\begin{equation*}
\Gamma(x,y)=\left\{
\begin{array}{l}
\cfrac{1}{2\pi}\ln(|x-y|) \text{ if } n=2,\\
-\cfrac{1}{4\pi}\cfrac{1}{|x-y|} \text{ if } n=3.\\
\end{array}
\right.
\end{equation*}
 The integral equations applying to direct problem will be obtained from a study of the classical single-
 and double-layer potentials. We begin to introduce  the following
 operators
\begin{equation*}
\begin{array}{lrll}
S_{\partial\Omega\partial\omega}:
 & u&\mapsto &S_{\partial\Omega\partial\omega} u(x):=
\displaystyle\int_{\partial\Omega}\Gamma(x,y) u(y) ~d\sigma(y);\\
S_{\partial\omega \partial\Omega}: &u&\mapsto& S_{\partial\omega
\partial\Omega}u(x):=
\displaystyle\int_{\partial \omega}\Gamma(x,y) u(y)~d\sigma(y)  ;\\
K_{\partial\Omega \partial\omega} : &u&\mapsto& K_{\partial\Omega
\partial\omega} u(x):= \displaystyle\int_{\partial
\Omega}\partial_\n\Gamma(x,y) u(y)~d\sigma(y)
~ ;\\
K_{\partial\omega\partial\Omega}: &u&\mapsto&
K_{\partial\omega\partial\Omega}u(x):= \displaystyle\int_{\partial
\omega}\partial_\n \Gamma(x,y) u(y) ~d\sigma(y)
\end{array}
\end{equation*}
Note that all these operators have a smooth kernel since the
boundaries $\partial\omega$ and $\partial\Omega$ are assumed to
have no common point. We  also denote
\begin{equation*}
\begin{array}{lrll}
S_{\Omega}:& u&\mapsto&S_{\Omega} u(x):=
\displaystyle\int_{\partial\Omega}\Gamma(x,y) u(y)
~d\sigma(y) ;\\
 K_{\Omega}:
&u&\mapsto& K_{\Omega}u(x):= \displaystyle\int_{\partial \Omega}\partial_\n \Gamma(x,y)u(y) ~d\sigma(y) ;
 \end{array}
\end{equation*}
\begin{equation*}
\begin{array}{lrll}
S_{\omega}: &u&\mapsto&S_{\omega}u(x):=
\displaystyle\int_{\partial\omega}\Gamma(x,y) u(y)
~d\sigma(y) ;\\
 K_{\omega} :
 &u &\mapsto&K_{\omega}u(x):= \displaystyle\int_{\partial \Omega}\partial_\n \Gamma(x,y)u(y) ~d\sigma(y) .
 \end{array}
\end{equation*}
%We refer to \cite{MazjaShraposhnikova} to precise the higher regularity properties of these operators. We quote the results we will need in the sequel.
%\begin{theorem}
%\label{resultat:mazjashaposhnikova}
%Let $l$ be a real number such that $l>(n+1)/2$. If $\partial \Omega$ belongs to $\sH^l$ (this means that locally $\partial\Omega$ is the graph of a $\sH^l$ function), then the following assertions hold
%\begin{itemize}
%\item the operator $\frac{1}{2} I +K_{\Omega}$ is an isomorphism of $\sH^l(\partial\Omega)$;
%\item the operator $S_{\Omega}$ maps isomorphically $\sH^{l-1}(\partial\Omega)$ onto $\sH^l(\partial\Omega)$.
%\end{itemize}
%\end{theorem}

We now obtain some systems of integral equations to compute the
state function and their shape derivatives. Since $v$ is harmonic
in $\Omega\setminus \overline{\omega}$ and for all $ x \in
\partial \Omega \cup \partial \omega $, it has  the classical
boundary representation:
\begin{equation}
\label{representationderive:Omega-omega}
\cfrac{1}{2} v(x) = \int_{\partial\Omega} \partial_\n\Gamma(x,y) v(y)- \int_{\partial\omega}
\partial_\n\Gamma(x,y) v(y)- \int_{\partial\Omega} \Gamma(x,y)  \partial_\n v(y)+ \int_{\partial\omega} \Gamma(x,y) \partial_\n v(y).
\end{equation}
Similarly since $v$ harmonic in $\omega$,  for all $ x \in
\partial \omega$ we can write
\begin{equation}
\label{representationderive:omega}
\cfrac{1}{2} v(x) =  \int_{\partial\omega}\partial_\n\Gamma(x,y) v(y)- \int_{\partial\omega} \Gamma(x,y) \partial_\n v(y).
\end{equation}
Let us denote by   $v_d$ the  solution of the boundary values
problem (D) in \eqref{probleme:transmission}. Let us show how to
compute their restrictions and also  their normal vector
derivatives on the boundaries. Incorporating the jump conditions,
a straightforward computation leads to the following  boundary
integral equations
\begin{equation}
\label{systemederivee:ud}
\begin{split}
\begin{bmatrix}
\cfrac{1}{2}I+\mu K_{ \omega} & \cfrac{\sigma_1}{\sigma_2+\sigma_1} S_{\partial\Omega \partial \omega}\\
\mu K_{\partial \omega\partial
\Omega}&\cfrac{\sigma_1}{\sigma_2+\sigma_1} S_{\Omega}
\end{bmatrix} &
\begin{bmatrix}
( \vd^+)_{|\partial\omega}\\
\\
\\
(\partial_\n \vd)_{|\partial\Omega}
\end{bmatrix}  \\=
\cfrac{1}{\sigma_1+\sigma_2}
&\begin{bmatrix}
\sigma_2\left(\cfrac{1}{2}I-K_{\omega}\right) &  S_{\omega}\\
\\
-\sigma_2 K_{\partial \omega\partial\Omega}& S_{\partial \omega \partial \Omega}
\end{bmatrix}
\begin{bmatrix}
\alpha
\\
\\
\\
\beta
\end{bmatrix}
+
\cfrac{\sigma_1}{\sigma_1+\sigma_2}
\begin{bmatrix}
 K_{\partial \Omega\partial\omega} f_1
\\
\\
\left(-\cfrac{1}{2}+K_{\Omega} \right) f_1
\end{bmatrix}
\end{split}
\end{equation}
where $\mu=[\sigma]/(\sigma_1+\sigma_2)$. Thanks to  \eqref{representationderive:omega}, the  quantity $(\partial_\n \vd)^+$ is then given by
\begin{equation*}
\label{partialnud2:omega}
 S_{\omega} (\partial_\n \vd)^+_{|\partial \omega} = \cfrac{\sigma_2}{\sigma_1} \left( -\cfrac{1}{2}I + K_{\omega} \right)\left(\vd^+(x)_{|\partial \omega}-\alpha\right)+ \cfrac{1}{\sigma_1} S_{\omega} \beta.
\end{equation*}
Concerning $\vn$, the solution of the Neumann problem (N) in
\eqref{probleme:transmission}, the same kind  of computations
gives
\begin{equation}
\begin{split}
\label{systemederivee:un}
\begin{bmatrix}
\cfrac{1}{2}I+\mu K_{ \omega} & -\cfrac{\sigma_1}{\sigma_2+\sigma_1} K_{\partial\Omega \partial \omega} \\
\mu K_{\partial \omega\partial
\Omega}&-\cfrac{\sigma_1}{\sigma_2+\sigma_1} \left(-\cfrac{1}{2}I+ K_{\Omega}\right)
\end{bmatrix}
&
\begin{bmatrix}
(\vn^+)_{|\partial\omega}\\
\\
\\
(\vn)_{|\partial\Omega}
\end{bmatrix} \\
= \cfrac{1}{\sigma_1+\sigma_2}&
\begin{bmatrix}
\sigma_2\left(\cfrac{1}{2}I-K_{\omega}\right) & S_{\omega}\\
\\
- \sigma_2 K_{\partial \omega\partial\Omega}&  S_{\partial \omega \partial \Omega}
\end{bmatrix}
\begin{bmatrix}
\alpha
\\
\\
\\
\beta
\end{bmatrix}
-\cfrac{\sigma_1}{\sigma_1+\sigma_2}
\begin{bmatrix}
 S_{\partial \Omega\partial\omega} g_1
\\
\\
S_{\Omega} g_1
\end{bmatrix}
\end{split}
\end{equation}
Finally, the computation of   $(\partial_\n \vn)^+_{|\partial \omega} $ is given by
\begin{equation*}%\label{partialnun2:omega}
 S_{\omega} (\partial_\n \vn)^+_{|\partial \omega} = \cfrac{\sigma_2}{\sigma_1} \left(-\cfrac{1}{2}I + K_{\omega} \right)\left(\vn^+(x)_{|\partial \omega}-\alpha\right)+
\cfrac{1}{\sigma_1} S_{\omega}\beta .
\end{equation*}
Concerning the well-posedness of  \eqref{systemederivee:ud}, we can state the following result.
\begin{theorem}
\label{theoreme:existence:systeme:integral}
The linear system of integral equation \eqref{systemederivee:ud}
has an unique solution in $\sH^{1/2}(\partial\omega) \times \sH^{-1/2}(\partial\Omega)$.
\end{theorem}

\begin{proofof}{Theorem \ref{theoreme:existence:systeme:integral}}
Let $A$ be the matricial operator defined on $\sH^{{1/ 2}}(\partial \omega)\times
\sH^{-1/2}(\partial \Omega)$ as
 \begin{equation}
A=
\begin{bmatrix}
\cfrac{1}{2}I+\mu K_{ \omega} & \cfrac{\sigma_1}{\sigma_2+\sigma_1} S_{\partial\Omega \partial \omega}\\
\mu K_{\partial \omega\partial
\Omega}& \cfrac{\sigma_1}{\sigma_1+\sigma_2}S_{\Omega}
\end{bmatrix}
\end{equation}
%Then, we can write equivalently the equations
%\begin{equation}
%\label{systeme:ud}
%A\begin{bmatrix}
%(\ud)_{|\partial\omega}\\
%\\
%\\
%(\partial_\n\ud)_{|\partial\Omega}
%\end{bmatrix} =
%\begin{bmatrix}
% \cfrac{\sigma_1}{\sigma_2+\sigma_1} K_{\partial \Omega\partial \omega}f\\
%\\
%\left(-\cfrac{1}{2}I+K_{ \Omega}\right)f
%\end{bmatrix}.
%\end{equation}
The main argument of the proof is based on the Fredholm alternative. In a first step, we  have to  show that the adjoint operator $A^*$ is injective. Since the boundaries are bounded, the adjoint operator $A^*$ defined on $\sH^{{-1/ 2}}(\partial \omega)\times
\sH^{1/2}(\partial \Omega)$
can be written under the form
\begin{equation}
A^*=
\begin{bmatrix}
\cfrac{1}{2}I+\mu K^*_{ \omega}
& \mu K^*_{\partial \Omega\partial
\omega}\\
\cfrac{\sigma_1}{\sigma_2+\sigma_1} S_{\partial\omega \partial
\Omega}&\cfrac{\sigma_1}{\sigma_2+\sigma_1} S_{\Omega}
\end{bmatrix}.
\end{equation}
Let $(u,v) \in \sH^{{-1/2}}(\partial \omega)\times \sH^{{1/2}}(\partial
\Omega)$ be in the kernel of  $A^*$. Consider the potential $W$ defined for each $x \in \R^d$ by
\begin{equation}
W(x)=\cfrac{\sigma_1}{\sigma_2+\sigma_1}\left( \int_{\partial\omega}
\Gamma(x,y) u(y) + \int_{\partial\Omega}\Gamma(x,y)v(y)
 \right).
\end{equation}
In a first step, we show that $W=0$. The function $W$ satisfies
$\Delta W=0$ in $\R^d\setminus(\partial{\omega} \cup \partial \Omega)$  by construction. We check that $W|_{\partial\Omega}=0$ from the equation corresponding to the second line
of $A^*$. By the properties of the single layer potential, $[W]=0$ on $\partial \omega$.
Furthermore, it holds $[\sigma{\partial_\n W}]=0$ on $\partial \omega$. Indeed, we can have (\cite{Kress})
\begin{equation*}
{\partial_\n W}^+=  \cfrac{\sigma_1}{\sigma_1+\sigma_2}\left(\left( \frac{1}{2}+
K^*_{\omega}\right)u+K^*_{\partial\Omega\partial\omega}v \right),
\end{equation*}
and
\begin{equation*}
{\partial_\n W}^-= \cfrac{\sigma_1}{\sigma_1+\sigma_2} \left(\left( - \frac{1}{2}+
K^*_{\omega}\right)v+K^*_{\partial\Omega\partial\omega}v \right),
\end{equation*}
hence,
\begin{equation*}
\sigma_1{\partial_\n W}^+-\sigma_2{\partial_\n W}^-=\sigma_{1}\left(
(\cfrac{1}{2}I+\mu K^*_{\omega})u
+\mu K^*_{\partial \Omega\partial
\omega}v \right).
\end{equation*}
This corresponds to the first line of $A^*(u,v)$. Then, $W$ solves the Laplace equation (\ref{equation:etat}) with
homogeneous Dirichlet boundary conditions. By the uniqueness of the
solution, we get $W=0$ in $\Omega$.

In a second step, we deduce that $u=v=0$. Since $W=0$ in
$\Omega$, we see that $[\partial_\n W]=0$ on $\partial\omega$. Since  $[\partial_\n W]=\sigma_1 u /(\sigma_1 +\sigma_2)$ on
$\partial\omega$ ,   we deduce $u=0$. From the second line of
$A^*(u,v)=0$, we see that  $S_{\Omega} v=0$ on $\partial \Omega$.
Since the single layer potential operator $S_\Omega  :
\sH^{-1/2}(\partial\Omega)\mapsto \sH^{1/2}(\partial\Omega)$ is
an isomorphism,  $v=0$ holds. The injectivity of $A^*$ is
proved.
Since $2A=I+C$ where $C$ is a compact operator, we conclude that $A$ has a continuous inverse thanks to the
Fredholm alternative.
\end{proofof}

 In a similar way, the problem \eqref{systemederivee:un} is well-posed under some additional assumptions.
 We define the adequate space
\begin{equation*}
\sH_\diamondsuit^{1/2}(\partial \Omega)=\left\{ \phi \in \sH^{1/2}(\partial \Omega)~:~ \int_{\partial \Omega} \phi =0\right\}.
\end{equation*}
We can state the following result.

\begin{theorem}
\label{existence:solution:second:systeme}
If we impose the normalizing condition
\begin{equation*}
\int_{\partial\Omega} v_n  =\int_{\partial\Omega}  f_1
\end{equation*}
then  there exists one unique couple
$((\vn)_{|\partial\omega},(\vn)_{|\partial\Omega}) \in \sH^{1/2}(\partial \omega)\times \sH^{1/2}_\diamondsuit(\partial \Omega)$ solution of \eqref{systemederivee:un} .
\end{theorem}

\begin{proofof}{Theorem \ref{existence:solution:second:systeme}}  Set
\begin{equation} B=
\begin{bmatrix}
\cfrac{1}{2}I+\mu K_{\omega} & - \cfrac{\sigma_1}{\sigma_2+\sigma_1} K_{\partial\Omega \partial \omega}\\ \mu K_{\partial \omega\partial
\Omega}&- \cfrac{\sigma_1}{\sigma_1+\sigma_2}\left(-\cfrac{1}{2}I+K_{\Omega}\right)
\end{bmatrix}
\end{equation}
the operator defined on $\sH^{{1/ 2}}(\partial \omega)\times
\sH_\diamondsuit^{1/2}(\partial \Omega)$. The adjoint $B^*$ can be written under the form
\begin{equation} B^*=
\begin{bmatrix}
\cfrac{1}{2}I+\mu K^*_{\omega} & \mu  K^*_{\partial\Omega \partial \omega}\\
-\cfrac{\sigma_1}{\sigma_1+\sigma_2}K^*_{\partial \omega\partial
\Omega}& -\cfrac{\sigma_1}{\sigma_1+\sigma_2}\left(-\cfrac{1}{2}I+K^*_{\Omega}\right).
\end{bmatrix}
\end{equation}
In a first step, we begin to show that $B^*$ is injective. Let $(u,v) \in \sH^{{1/2}}(\partial \omega)\times \sH^{{1/2}}(\partial
\Omega)$ be in the kernel of  $B^*$. We introduce  the potential
\begin{equation*}
Z(x)=- \frac{\sigma_1}{\sigma_1+\sigma_2}\left(\int_{\partial
\omega}\Gamma(x,y)u(y) +\int_{\partial
\Omega}\Gamma(x,y)v(y) \right),~x\in {\mathbb{R}^d}.
\end{equation*}
We can see that $Z$ is a harmonic function in
${\mathbb{R}}^d\backslash (\partial\omega\cup\partial\Omega)$,
satisfying $\partial_\n Z~~ =0$ on $\partial \Omega$. By the properties of the single layer potential, $[Z]=0$ Furthermore, a straightforward calculation shows that  $[\sigma \partial_\n Z]=0$ on $\partial\omega$. Hence, $Z$ solves the boundary value problem
\begin{eqnarray*}
-\div{\sigma \nabla Z} &=& 0 \text{ in } \Omega,\\
\partial_n Z &=& 0 \text{ on } \partial \Omega.
\end{eqnarray*}
The function is therefore constant in $\Omega$.
Writing $[\partial_\n Z]=0$ on $\partial\omega$, we get
easily $u=0$ and then $(- \frac{1}{2}+K^*_\Omega)v=0$. Since the
operator $\lambda I-K^*_\Omega$ is one to one on
$\sH_\diamondsuit^{1/2}(\partial \Omega)$, we deduce that $v=0$.
 We conclude the proof thanks to the Fredholm
alternative.
\end{proofof}

\textbf{Tangential regularity results.}
Let us consider now the particular case where both $\alpha$ and
$\beta$ are the zero function and $(f_1,g_1)=(f,g)$ where $f$ and
$g$ are respectively the Dirichlet and Neumann boundary data. To
recover the tangential regularity of the solution $u$ along
$\partial\omega$, we look at  the first line of \eqref{systemederivee:ud} to deduce that 
\begin{equation}
 \label{equation:revue}
\left[ \cfrac{1}{2}I+\mu K_{ \omega} \right]
(\ud)_{|\partial\omega} =  - \cfrac{\sigma_1}{\sigma_2+\sigma_1}
S_{\partial\Omega \partial \omega}\partial_n \ud|_{\partial\Omega}
+ \cfrac{\sigma_1}{\sigma_1+\sigma_2} K_{\partial
\Omega\partial\omega} f;
\end{equation}
and
\begin{equation}
S_{\omega} (\partial_\n \ud)^+_{|\partial \omega} =
\cfrac{\sigma_2}{\sigma_1} \left( -\cfrac{1}{2}I + K_{\omega}
\right)\ud^+(x)_{|\partial \omega}
\end{equation}
It is easy to deduce that $(\ud)_{|\partial\omega} \in
\mathcal{C}^{3,\alpha}(\partial\omega)$. Indeed, from
 \eqref{equation:revue} that we consider  as an equation in
$(\ud)_{|\partial\omega}$ with data $f$ and
$(\partial_\n\ud)_{|\partial\Omega}=g$, we see  that
$(f,(\partial_\n\vd)_{|\partial\Omega})$ belongs to
$\sH^{1/2}(\partial\Omega)\times \sH^{-1/2}(\partial\Omega)$, thanks to Theorem
\ref{theoreme:existence:systeme:integral}.

In order to give a sense to the jump conditions arising in \eqref{derivee:etatD},\eqref{derivee:etatN},\eqref{derivee:seconde:etatD}, we need to work in space of functions of higher regularity. We choose the framework of H\"{o}lder spaces. We quote \cite{KirschM2AS} to precise the behavior of the layer potentials on these spaces.
 
\begin{theorem}(Kirsch \cite{KirschM2AS})
\begin{enumerate}
\item
If $\partial \omega$ is of class $C^{2,\alpha}, ~0<\alpha<1$
then the operators $S_\omega$ and $K_\omega$ map
$C^{\beta}(\partial \omega)$ continuously into $C^{1,\beta}$ for
all $0<\beta\le \alpha$. \item Let $k \in \mathbb{N}$ with $k \neq
0$. If $\partial \omega$ is of class $C^{k+1,\alpha}$ with
$0<\alpha<1$, then the operators $S_\omega$ and $K_\omega$ map
$C^{k,\beta}(\partial \omega)$ continuously into
$C^{k+1,\beta}(\partial \omega)$ for all $0<\beta\le \alpha$.
\item Let $k$ be an integer. If $\partial \omega$ is of class
$C^{k+2,\alpha}$, then $K^*_\omega$ maps $C^{k,\beta}$
continuously into $C^{k+1,\beta}(\partial\omega)$ for all
$0<\beta\le \alpha$.
\end{enumerate}
\end{theorem}

We go back to the proof. Since the two boundaries have no intersection point and since $\partial \omega$ is of class $C^{4,\alpha}$, it  follows that the
right hand side of the former equation is of class
$\mathcal{C}^{3,\alpha}(\partial\omega)$. We then conclude  the
solution of \eqref{equation:revue} will be of class
$\mathcal{C}^{3,\alpha}$ since  the operator $1/2I+\mu K_{ \omega}
$ is an isomorphism from $\mathcal{C}^{3,\alpha}(\partial\omega)$
into itself. With the same arguments, we show straightforwardly
that $(\partial_\n \un)^+_{|\partial \omega} \in \mathcal{C}^{2,\alpha}$.
\par
\noindent
\textbf{About the regularity of the jumps of the second derivative.}~
The equations giving the jump conditions $[\ud']$ and $[\partial
\ud']$ show  obviously that $[\ud']$ and $[\partial_n \ud']$
belong respectively to $C^{2,\alpha}(\partial \omega)$ and
$C^{1,\alpha}(\partial \omega)$. Hence, it comes straightforwardly
that  $[\ud''] \in \mathcal{C}^{1,\alpha}$. With the same arguments, we show
that $[\partial_n \ud''] \in \mathcal{C}^{0,\alpha}$(see
\cite{SokolowskiZolesio} for more details) and then  that all the
formal computations to get the equations describing the second
derivative have a sense.

\begin{remark}
In a view of a numerical discretization of the state equation, one has to emphasize that the choice of a finite elements method seems inappropriate: one should extract tangential derivative of high order on the interface $\partial\omega$. The obtained numerical accuracy is not sufficient to incorporate the results in an optimization scheme. On the converse, the systems  of boundary integral equations \eqref {systemederivee:ud} and \eqref{systemederivee:un} are well-suited for this kind of  computation. Nevertheless, a discussion of adapted schemes should be precise  and is out of the scope of this manuscript.
\end{remark}

\subsection{Case of Neumann boundary conditions.}

Since the admissible deformation fields have a support with no
intersection points with the outer boundary, it is a
straightforward application of the preceding computations to show
that $\un$ solution to
\eqref{equation:etat}-\eqref{equation:potentiel:condition:limite:neumann}
is twice differentiable with respect to the shape. Furthermore,
its second order derivative $\un''$ belongs to
$\sH^1(\Omega\setminus\overline{\omega})\cup\sH^1(\omega)$  and
solves
\begin{equation}
\label{derivee:seconde:etatN}
\left\{
\begin{array}{rcl}
\Delta \un''  & =  & 0 \text{ in }\omega\setminus\overline{\omega} \text{ and  in } \omega,   \\
  \left[\un''\right]& =&  \left( h_{1,n} h_{2,n}H -\ct{{\h}_1}.D\n\, \ct{{\h}_2} \right) [\partial_{\n} \un] - \left( h_{1,n}[\partial_{\n} (\un)'_2] + h_{2,n}[\partial_{\n} (\un)'_1] \right) \\
  && ~~~~~~+ \left(\ct{{\h}_1}.\nabla h_{2,n} +\ct{{\h}_2}.\nabla h_{1,n} \right) [\partial_{\n} \un] \text{ on } \partial \omega, \\
  \left[\sigma \partial_{n} \un'' \right]&  = & \divt{h_{2,n}\left[\sigma \nablat (\un)'_1 \right] + h_{1,n}\left[\sigma \nablat (\un)'_2 \right] +\ct{h_1}. D\n .\ct{h_2})[\sigma \nablat \un]}\\
 && -  \divt{(\ct{{\h}_1}.\nablat h_{2,n} +  \nablat h_{1,n}.\ct{h_2} +  h_{2,n} h_{1,n}  (2  D\n-H I))  \left[\sigma \nablat \un \right]   } \text{ on } \partial \omega,\\
   \partial_n\un'' &=& 0 \text{ on } \partial\Omega;
\end{array}
\right.
\end{equation}
where we use the notations of Theorem \ref{theoreme:derivee:seconde:etatD}.
%%%%%%%%%%%%%%%%%%%%%%%%%%%%%
\section{Second order derivatives for the criterion.}
\label{section:derivee:seconde:critere}
%%%%%%%%%%%%%%%%%%%%%%%%%%

%%%%%%%%%%%%%%%%%%%%%%%%%%%%%
\subsection{Proof of Theorem \ref{derivee:seconde:kohn:vogelius}.}
%%%%%%%%%%%%%%%%%%%%%%%%%%%%%

The differentiability of the objective is a direct application of Theorem \ref{theoreme:derivee:seconde:etatD}.
The computation we make here is based on the relation
\begin{equation}\label{relation:derivee:seconde:fonctionnelle}
D^2\KV(\omega)(\h_1,\h_2)=D \left(D\KV(w)\h_1\right)\h_2 -D\KV(w) D{\h_1}\h_2.
\end{equation}
To obtain \eqref{second:derivative}, we  compute in a first step the shape
gradient in the direction $\h_1$. Then, in a second step, we  differentiate the obtained expression in the direction of $\h_2$. In the sequel, we adopt the notation $v
=\ud-\un$ to obtain concise expressions.
\begin{eqnarray*}
D\KV(\omega)\h_1 & =&
\sigma_1 \int_{\Omega\setminus\overline{\omega}} \div{ |\nabla v|^2 \h_1}+2 \nabla v.\nabla v'_1 +
 \sigma_2 \int_{\omega} \div{|\nabla v|^2 \h_1}+2 \nabla v.\nabla v'_1\\
&=& \sigma_1 (A_1+2 B_1) +\sigma_2 (A_2+ 2 B_2),
\end{eqnarray*}
where
\begin{eqnarray*}
A_1&=&\int_{\Omega\setminus\overline{\omega}} \div{ |\nabla v|^2 \h_1}~~~~~~B_1=\int_{\Omega\setminus\overline{\omega}} \nabla v.\nabla v'_1\\
A_2&=&\int_{\omega} \div{|\nabla v|^2 \h_1}~~~~~~~~B_2=\int_{\omega} \nabla v.\nabla v'_1.
\end{eqnarray*}
Now, we use the classical formulae to differentiate a domain integral to get
\begin{eqnarray*}
DA_1(\omega)\h_2&= & \int_{\Omega\setminus\overline{\omega}}
\div{\div{ |\nabla v|^2 \h_1} \h_2}  +2\div{ \nabla v.\nabla v'_2 ~\h_1},\\
 & =& - \int_{\partial\omega} \div{ |\nabla v^+|^2 \h_1} h_{2,n}+ 2  \nabla v^+.\nabla (v'_2)^+ ~h_{1,n};\\
DA_2(\omega)\h_2 &= & \int_{\partial \omega}
\div{ |\nabla v^-|^2 \h_1} h_{2,n}+ 2  \nabla v^-.\nabla (v'_2)^- ~h_{1,n}.
\end{eqnarray*}
The terms $DB_i,~i=1,2$ require more precisions. First, we write
\begin{eqnarray*}
DB_1(\omega)\h_2 & = & \int_{\Omega\setminus\overline{\omega}}  \div{ \nabla v .\nabla v'_1 ) \h_2} + \nabla v'_1 . \nabla v'_2  + \nabla v. \nabla (v'_1)'_2, \\
&=& -\int_{\partial\omega} \nabla v^+ . \nabla (v'_1)^+ h_{2,n} + \partial_\n v^+ ((v'_1)^+)'_2 + \cfrac{1}{2} \left(\partial_\n (v'_1)^+ (v'_2)^+  +   \partial_\n (v'_2)^+ (v'_1)^+\right) \\
&&- \int_{\partial\Omega} \partial_\n v ((\un)'_1)'_2 + \cfrac{1}{2} \left(\partial_\n (\ud)'_1 (\un)'_2 + \partial_\n (\ud)'_2 (\un)'_1 \right)
\end{eqnarray*}
Note that we used the Green formula twice to keep the symmetry in
$\h_1$ and $\h_2$. We also use the fact that the derivatives
$(\ud)'_i$ are harmonic in $\Omega\setminus\overline{\omega}$ to
transform the boundary integral on the exterior boundary into an
integral on the moving boundary. We obtain
\begin{eqnarray*}
DB_1(\omega)\h_2 & = &- \int_{\partial\omega} \nabla v^+. \nabla (v'_1)^+ \h_{2,n} +  \partial_\n v^+ (((\ud)'_1)'_2)^+ - v \partial_\n (((\un)'_1)'_2)^+\\
&-&\int_{\partial\omega}\cfrac{1}{2} \left( \partial_\n (v'_1)^+ ((\ud)'_2)^+ + \partial_\n (v'_2)^+ ((\ud)'_1)^+
    - \partial_\n (\un)'_1)^+ (v'_2)^+ - \partial_\n ((\un)'_2)^+ (v'_1)^+ \right).
\end{eqnarray*}
By the same methods, we get
\begin{eqnarray*}
DB_2(\omega)\h_2 & = & \int_{\partial\omega} \nabla v^-. \nabla (v'_1)^- h_{2,n} + \partial_\n v^- ((v'_1)'_2)^- + \cfrac{1}{2} \left( \partial_\n (v'_1)^- (v'_2)^- + \partial_\n (v'_2)^- (v'_1)^- \right).
\end{eqnarray*}
We regroup the different terms and  after some  straightforward computations, we obtain:
\begin{eqnarray*}
D \left(D\KV(\omega)\h_1\right)(\omega) \h_2&= &
- \int_{\partial \omega} \div{\left[\sigma | \nabla v |^2 \h_1 \right]} +
2 \left[\sigma \nabla v . \left( h_{1,n} \nabla v'_2 + h_{2,n} \nabla v'_1 \right)\right]\\
%&&- \int_{\partial \omega} \left[\sigma \left( \partial_\n  (\ud)'_1 (\ud)'_2 + \partial_\n  (\ud)'_2 (\ud)'_1  +   \partial_\n  (\un)'_1 (\un)'_2  +\partial_\n  (\un)'_2 (\un)'_1\right)\right] \\
%&&+\int_{\partial \omega} \left[ \sigma \left( \partial_\n  (\un)'_1 (\ud)'_2  + \partial_\n  (\un)'_2 (\ud)'_1 + \partial_\n  (\un)'_2 (\ud)'_1 + \partial_\n  (\un)'_1 (\ud)'_2 \right)\right] \\
&&- \int_{\partial \omega} \left[\sigma \left( (\ud)'_2 \partial_\n  v'_1 + (\ud)'_1   \partial_\n  v'_2 - \partial_\n  (\un)'_2 v'_1- \partial_\n  (\un)'_1 v'_2  \right)\right] \\
&& +2  \int_{\partial \omega} v \left[ \sigma \partial_\n ((\un)'_1)'_2 \right] - \sigma_1 \partial_\n v^+ \left[((\ud)'_1)'_2)\right] .
\end{eqnarray*}
 In order to compute $D^2\KV(\omega)(\h_1,\h_2)$,  the first order derivative of the Kohn-Vogelius objective is needed. It can be written  as follows:
\begin{equation*}
D\KV(w)\h=- \int_{\partial \omega} \left[\sigma| \nabla v |^2  \right] h_n+2\int_{\partial \omega} v\left[ \sigma \partial_\n (\un)'\right] -\sigma_1 \partial_\n v^+ \left[(\ud)'\right].
\end{equation*}
Gathering
\eqref{relation:derivee:seconde:fonctionnelle},\eqref{relation:potentiel:derivee:seconde}
and \eqref{relation:flux:derivee:seconde}, we write the second
derivative of the Kohn-Vogelius criterion as:
 \begin{eqnarray*}
D^2\KV(\omega)(\h_1,\h_2)&= &- \int_{\partial \omega} \div{\left[\sigma | \nabla v |^2 \h_1 \right]} - \left[\sigma | \nabla v |^2 \right]( D{\h_1}\h_2).\n\\
&&- \int_{\partial \omega} \left[\sigma \left( (\ud)'_2 \partial_\n  v'_1 + (\ud)'_1   \partial_\n  v'_2 - \partial_\n  (\un)'_2 v'_1- \partial_\n  (\un)'_1 v'_2  \right)\right] \\
&& +2  \int_{\partial \omega}  \left[\sigma \nabla v . \left( h_{1,n} \nabla v'_2 + h_{2,n} \nabla v'_1 \right)\right]+v \left[ \sigma \partial_\n (\un)''_{1,2} \right] - \sigma_1 \partial_\n v^+ \left[(\ud)''_{1,2}\right] .
\end{eqnarray*}
Let us give a more simplified version for  the first term.  We
decompose the field $\h_2$ into normal vector and tangential parts
and we  use \eqref{Jacobien:suivant:tangente}. After some
elementary computations, we obtain
\begin{equation*}
\begin{split}
- \int_{\partial \omega}& \div{\left[\sigma | \nabla v |^2 \h_1 \right]} -\left[\sigma | \nabla v |^2 \right]( D{\h_1}\h_2).\n\\&=\int_{\partial \omega} \left[\sigma | \nabla v |^2 \right] \left( \ct{h_1}.\nablat h_{2,n}+ \ct{h_2}.\nablat h_{1,n}-\ct{h_2}. D\n\, \ct{h_1}\right)-\int_{\partial \omega} \partial_\n \left(\left[\sigma | \nabla v |^2  \right]\right) h_{1,n} h_{2,n}.
\end{split}
\end{equation*}
Finally, the second order derivative of the Kohn-Vogelius objective becomes:
\begin{equation}
\begin{split}\label{derivee:seconde:fonctionnelle:frontiere}
D^2\KV(\omega)(\h_1,\h_2)&= \int_{\partial \omega} \left[\sigma | \nabla v |^2 \right] \left( \ct{h_1}.\nablat h_{2,n}+ \ct{h_2}.\nablat h_{1,n}-\ct{h_2}. D\n\, \ct{h_1}\right)\\
&-\int_{\partial \omega} \partial_\n \left(\left[\sigma | \nabla v |^2  \right]\right) h_{1,n} h_{2,n}+2 \left[\sigma \nabla v . \left( h_{1,n} \nabla v'_2 + h_{2,n} \nabla v'_1 \right)\right]\\
&- \int_{\partial \omega} \left[\sigma \left( (\ud)'_2 \partial_\n  v'_1 + (\ud)'_1   \partial_\n  v'_2 - \partial_\n  (\un)'_2 v'_1- \partial_\n  (\un)'_1 v'_2  \right)\right] \\
&+2  \int_{\partial \omega} v \left[ \sigma \partial_\n (\un)''_{1,2} \right] - \sigma_1 \partial_\n v^+ \left[(\ud)''_{1,2}\right] .
\end{split}
\end{equation}

%%%%%%%%%%%%%%%%%%%%%%%%%%%%%%%%%%%%
\subsection{Analysis of stability. Proof of Theorem \ref{compactness:D2KV}}
%%%%%%%%%%%%%%%%%%%%%%%%%%%%%%%%%%%%

Now, we specify the domain $\omega$ that is assumed to be a critical shape for $\KV$. Moreover, we assume that the additional condition $\ud=\un$ holds. To emphasize that we deal with such a special domain, we will denote it  $\omega^*$. The assumptions mean that the measurements are compatible and that $\omega^*$ is a global minimum of the criterion. From the necessary condition of order two at a minimum, the shape Hessian is positive at such a point.
 
Let us notice that only the normal component of $\h$ appears. Let us also emphasize that there is no hope to get
$\h=0$ from the structure theorem for second order shape
derivative (\cite{HenrotPierre}). The deformation
field $\h$ appears in $D^2\KV(\omega^*)(\h,\h)$ only thought its
normal component $h_n$ since $\omega^*$ is a critical point for
$\KV$. This remark explains why we consider in the statement of
Theorem \ref{compactness:D2KV} the  scalar Sobolev space
corresponding to the normal components of the deformation field.

We now prove Theorem \ref{compactness:D2KV}. From \eqref{derivee:seconde:fonctionnelle:frontiere}, we deduce
\begin{equation}
\begin{array}{llll}
 D\KV^2(\omega^*)[h,h]&=- 2\displaystyle\int_{\partial \omega^*} \left[\sigma \left( \ud' \partial_\n  v'- \partial_\n  \un'  v' \right)\right]\\
&= 2\left[\sigma \right]\displaystyle\int_{\partial \omega^*} \left((\ud'^+- \un'^+) \divt{h_n \nablat \ud}-\cfrac{\sigma_1}{\sigma_2}\partial_\n \ud^+ h_n \partial_\n(\ud'-\un')^+  \right)\\
&= 2\left[\sigma \right]\left(\Big\langle \ud'^+- \un'^+,\divt{h_n \nablat \ud}\Big\rangle-\cfrac{\sigma_1}{\sigma_2}\Big \langle \partial_\n \ud h_n ,\partial_\n\left(\ud'-\un'\right)^+ \Big\rangle\right).
\end{array}
\end{equation}
 where $\langle ,\rangle$ denotes the duality between $H^{1/2}(\partial\omega^*)\times H^{-1/2}(\partial\omega^*)$ .
Let us introduce the operators
\begin{equation*}
\begin{array}{rclrcl}
 T_1:H^{1/2}(\partial \omega^*) &  \rightarrow  &  H^{-1/2}(\partial \omega^*)   &M_1 :H^{1/2}(\partial \omega^*) &\rightarrow &  H^{1/2}(\partial \omega^*) \\
  \h& \mapsto &\divt{h_n \nablat \ud} &  \h& \mapsto  & \ud'^+- \un'^+\\
 T_2:H^{1/2}(\partial \omega^*) &  \rightarrow  &  H^{1/2}(\partial \omega^*)   &M_2 :H^{1/2}(\partial \omega^*)&\rightarrow &  H^{-1/2}(\partial \omega^*) \\
  \h& \mapsto & h_n \partial_\n \ud^+  &  \h& \mapsto  & \partial_\n\left(\ud'^+- \un'^+\right)\\
\end{array}
\end{equation*}
The Hessian can then be written under the form :
\begin{equation*}
D^2\KV(\omega^*)(\h,\h) = 2\left[\sigma\right] \left(\Big\langle M_1(\h),T_1(\h)\Big\rangle-\cfrac{\sigma_1}{\sigma_2}\Big\langle T_2(\h), M_2(\h)\Big\rangle \right).
\end{equation*}
 From the classical results  of Maz'ya and Shaposhnikova on multipliers (\cite{VGMaz}, \cite{Triebel}), we  get easily that $T_1$ and $T_2$
 are continuous operators.  In
fact, the compactness of the Hessian is a consequence of the fact
that both operators $M_1$ and $M_2$ are compact. We use a
regularity argument : we remark that $M_1$ is the composition of
the operators:
\begin{equation*}
\begin{array}{rclcrcl}
R_1:  H^{1/2}(\partial \omega^*)   &  \rightarrow  &  H^{1/2}_{\diamond}(\partial \Omega)  &  \text{ and}  & R_2:  H^{1/2}_{\diamond}(\partial \Omega) &\rightarrow & H^{1/2}(\partial \omega^*)  \\
  \h& \mapsto &   -\un' & &\phi& \mapsto  &   \psi\\
\end{array}
\end{equation*}
where $\psi$ is the trace on $\partial\omega^*$ of $\Psi$ solution of
\begin{equation*}
\left\{
\begin{array}{rcl}
- \Delta \Psi  & =  & 0 \text{ in }\Omega\setminus\overline{\omega^*} \text{ and  in } \omega^*,  ,   \\
 \left[ \Psi\right]&=&0 \text{ on }\partial \omega^*,\\
  \left[\sigma\partial_\n \Psi\right] & =  & 0  \text{ on } \partial \omega^*,\\
   \Psi &=& \phi \text{ on } \partial\Omega.
\end{array}
\right.
\end{equation*}
While $R_1$ is a continuous operator, we prove that $R_2$ is compact.
Let us express $u_{|\partial \omega^*}= \psi$. We use the integral formula of $u$ to obtain:
\begin{equation*}
\begin{bmatrix}
\cfrac{1}{2}I+\mu K_{ \omega^*} & \cfrac{\sigma_1}{\sigma_2+\sigma_1} S_{\partial\Omega \partial \omega^*}\\
\mu K_{\partial \omega^*\partial
\Omega}&\cfrac{\sigma_1}{\sigma_2+\sigma_1} S_{\Omega}
\end{bmatrix}
\begin{bmatrix}
( u)_{|\partial\omega^*}\\
\\
\\
(\partial_\n u)_{|\partial\Omega}
\end{bmatrix} =
\cfrac{\sigma_1}{\sigma_1+\sigma_2}
\begin{bmatrix}
 K_{\partial \Omega\partial\omega^*} \phi
\\
\\
\left(-\cfrac{1}{2}+K_{\Omega} \right) \phi
\end{bmatrix}
\end{equation*}
The matricial operator arising in this equation appeared also
in $\eqref{systemederivee:ud}$. It  has a continuous inverse thanks to  Theorem
\ref{theoreme:existence:systeme:integral}. Let us express $u|_{\partial \omega^*}=\psi$:
\begin{equation}
\label{estimee:hessienne}
\left[(\cfrac{1}{2}I+\mu K_{\omega^*})- \mu S_{\partial \Omega\partial \omega^*}S_{\Omega}^{-1} K_{\partial \omega^*\partial \Omega}\right]\psi=\cfrac{\sigma_1}{\sigma_1+\sigma_2}\left[K_{\partial \Omega\partial \omega^*}-S_{\partial \Omega\partial \omega^*}S_{\Omega}^{-1}(\cfrac{1}{2}I-K_{\Omega})\right]\phi.
\end{equation}
Since  the operators $K_{\partial \Omega \partial \omega^*}$ and $S_{\partial \Omega\partial \omega^*}$ are compact, the operator  $R_2$ is compact, hence $M_1$  is compact. The proof of compactness of  $M_2$ is similar.  Let us mention that a similar strategy of proof can be found in \cite{EpplerHarbrecht}.

The natural question is then to quantify how is this optimization problem ill-posed. This question is directly in related to the rate at which the singular values of the Hessian operator are decreasing.
Equation \eqref{estimee:hessienne} shows that this rate is the one of the operators
$K_{\partial \Omega \partial \omega^*}$ and $S_{\partial \Omega\partial \omega^*}$.
Now, since for every $u\in \sH^{1/2}(\partial\Omega)$,
the functions $K_{\partial \Omega \partial \omega^*}u$ and
$S_{\partial \Omega\partial \omega^*}u$ are harmonic outside of $\partial\Omega$ and
therefore in $\Omega$, their restrictions on $\partial\omega^*$ are as smooth as
$\partial\omega^*$. We conclude that if $\partial\omega^*$ is $\mathcal{C}^{\infty}$ then the restriction belongs to
each $\sH^s(\partial\omega^*)$ for $s>1/2$  then that if $\lambda_n$ denotes the $n^{th}$ eigenvalue
of $D^2\KV(\omega^*)$, then $\lambda_n = o(n^{-s})$ for all $s>0$.

%\section{Appendix : Regularity results}

% Since we know that the injectivity is true in the space $L^2(\partial \omega)$, we
% 3),
% since $\omega$ and $\omega^c$ are connected,  Zaremba's principle enables us to show the injectivity and then the invertibility of $T_i,~i=1,2$.

%%%%%%%%%%%%%%%%%%%%%%%%%%%%%%%%%%%%%%%%%%%%%

\end{document}